\documentclass{article}
\usepackage{amsmath}  
\usepackage{amssymb}  
\usepackage{amsfonts} 
\usepackage{amsthm}   
\usepackage{bm}       
\usepackage{mathtools} 
\usepackage{cancel}
\usepackage{comment}
\usepackage{hyperref}
\usepackage[dvipsnames]{xcolor}
\usepackage{enumitem}
\usepackage{booktabs}
\usepackage{multirow}
\usepackage{adjustbox}
\usepackage{tikz}
\usepackage{amsmath}
\usetikzlibrary{patterns, decorations.pathreplacing}
\usepackage{array}
\usepackage{makecell}
\usepackage{authblk} 
\usepackage{abstract}
\usepackage{subcaption}
\usepackage{natbib}

\usepackage[a4paper,margin=1in]{geometry}

\newcommand{\E}{\mathbb{E}}
\newcommand{\V}{\mathbb{V}}
\newcommand{\C}{\mathbb{C}\mathrm{ov}}

\newcommand{\tas}{\xrightarrow{a.s.}}
\newcommand{\Prob}{\mathbb{P}}

\newtheorem{theorem}{Theorem}
\newtheorem{lemma}{Lemma}

\newtheorem{proposition}{Proposition}
\newtheorem{assumption}{Assumption}

\theoremstyle{definition} 

\newtheorem{remark}{Remark}

\makeatletter
\def\@biblabel#1{}
\makeatother

\newcommand{\tl}{\xrightarrow{\hspace{0.1cm}\mathcal{L}\hspace{0.1cm}}}
\newcommand{\tp}{\xrightarrow{\hspace{0.1cm}\Prob\hspace{0.1cm}}}

\date{}

\title{Testing the equality of estimable parameters across many populations}

\begin{document}
\author[]{M. Romero-Madroñal\thanks{Corresponding author. Email address: \texttt{mrmadronal@us.es}}}
\author[]{M. Remedios Sillero-Denamiel}
\author[]{M. Dolores Jiménez-Gamero}

\affil[]{Department of Statistics and Operations Research, Universidad de Sevilla, Spain \\ Instituto de Matem\'aticas de la Universidad de Sevilla, Spain}

\maketitle

\begin{abstract}
\normalsize
The comparison of a parameter in $k$ populations is a classical problem in statistics. Testing for the equality of means or variances are typical examples. Most procedures designed to deal with this problem assume that $k$ is fixed and that samples with increasing sample sizes are available from each population. This paper introduces and studies a test for the comparison of an estimable parameter across $k$ populations, when $k$ is large and the sample sizes from each population are small when compared with $k$. The proposed test statistic is asymptotically distribution-free under the null hypothesis of parameter homogeneity, enabling asymptotically exact inference without parametric assumptions. Additionally, the behaviour of the proposal is studied under alternatives. Simulations are conducted to evaluate its finite-sample performance, and a linear bootstrap method is implemented to improve its behaviour for small $k$. Finally, an application to a real dataset is presented.


\end{abstract}

\section{Introduction}


Testing the equality of a parameter in a given number of populations from independent samples coming from them, is a central topic in statistics, with early methods including one-way ANOVA for comparing means \citep{Fisher1925} and Bartlett’s test for comparing variances \citep{Bartlett1937}. Given the sensitivity of these procedures to parametric assumptions, nonparametric alternatives such as the Kruskal–Wallis' test \citep{Kruskal1952} and Levene’s test \citep{Levene1960} were later developed.

These classical methods implicitly assume that the number of populations $k$ is fixed and small relative to the within-group sample sizes $n_1,\ldots,n_k$, which are usually assumed to be large. In many contemporary applications, such as genomics (see e.g. \cite{ZhanHart2014}) or social sciences (see e.g. \cite{AlbaFernandez2024}), researchers frequently face problems where $k$ is large and larger than the sample sizes from each population. In an era of massive data collection, comparing fundamental parameters, such as measures of location, dispersion, or association, across many groups remains an interesting problem.

The behavior of some classical tests in the large-$k$, small-$n$ regime has been studied in the literature for some specific testing problems. Below, we cite some papers dealing with the comparison of means. \cite{BoosBrownie1995} analyzed the asymptotic properties of ANOVA and rank-based tests as $k \to \infty$, assuming homoscedasticity. Building on this work, \cite{AkritasArnold2000} derived the asymptotic distribution of the ANOVA $F$-statistic in both one-way and two-way layout models. Subsequently, \cite{Bathke2002} extended the analysis to the balanced multi-factor case, demonstrating that the asymptotic properties observed in simpler layouts can be generalized to more complex factorial designs. Later, \cite{AkritasPapadatos2004} introduced nonparametric procedures for one-way layouts that accommodate heteroscedasticity, thereby broadening the applicability of multi-sample inference in large-$k$ settings. Collectively, these studies show that classical $F$-statistics and their rank-based analogues, under the assumption of homoscedasticity, converge in distribution to a normal law as $k \to \infty$, deviating from the classical $F$ distribution. This shift (which not only happens in the tests for the means) highlights the  necessity of either adapting existing tests or developing novel statistical methods specially designed for the large-$k$, small-$n$ context.


Motivated by these challenges, some recent research is focused on multi-sample inference with large $k$ and small sample sizes for testing other hypothesis rather than comparing means.
Specifically, goodness-of-fit tests have been designed for this regime, with recent contributions addressing Poissonity, normality, and proportion tests in \cite{JimenezGamero2024a}, \cite{JimenezGamero2024b}, \cite{AlbaFernandez2024}, respectively; approaches for the $k$-sample problem have also been developed by \cite{ZhanHart2014}, \cite{Kim2021}, and \cite{JimenezGamero2022}; a test for homocedasticity has been proposed in \cite{JimenezGamero2025}. 

In particular, in their practical application, \cite{ZhanHart2014} rejected the null hypothesis of equal densities, suggesting the presence of heterogeneity across the populations. This naturally raises the question of whether the observed differences are driven by specific parameters such as the mean, variance, or other distributional characteristics. Investigating these aspects could provide a deeper understanding of the nature of the heterogeneity detected by the test.


While much of the existing work has focused on comparing means, variances, or entire distributions, it is natural to consider broader comparisons involving parameters such as measures of location, scale, dispersion, association, or other functionals expressible as expectations of kernels. In this way, the comparison of means and variances can be extended to other population parameters, such as the Gini mean difference and Spearman’s rho. In this direction, the framework developed in this paper provides a unified approach to compare such parameters across many populations without relying on specific structural properties of the kernel.


In this work, assuming independent samples from each population, the problem of testing the null hypothesis \( H_0: \theta_1 = \dots = \theta_k \) is considered, where each \( \theta_i \in \mathbb{R} \) is defined as the expectation of a kernel \( h \) of fixed degree \( m \geq 1  \). A test statistic \( T_k \) is proposed, which is an unbiased estimator of the variance of  the parameters $D_k = (1/k) \sum_{i=1}^k (\theta_i - \overline{\theta})^2$.
Notice that $H_0$ is equivalent to $D_k=0$. Under some mild assumptions, \( T_k \) is consistent and asymptotically normal, where by asymptotic we mean as $k \to \infty$. The asymptotic distribution of \( T_k \) is derived by decomposing it  into a linear component and a remainder. The variance of the reminder is shown to be 
is asymptotically negligible relative to that of the linear part. The practical use of $T_k$ as a test statistic, requires the estimation of its asymptotic variance, which is the variance of its linear component. A ratio-consistent estimator of the variance under the null hypothesis is devised. This way, it is obtained a test statistic which is asymptotically free distributed under the null.
The power of the proposed test is analyzed for a broad class of alternatives, establishing its consistency under explicit conditions on the rate at which alternatives converge to the null and the growth of the sample sizes. 

To improve the 
performance of the proposed procedure for moderate and small values of $k$  we additionally implement a bootstrap procedure. The behavior of the asymptotic approximation and the bootstrap approximation are assessed through  an extensive simulation study focusing on the Gini mean difference and Spearman’s rho.

One of the assumed assumptions to derive the previous results is that the sample sizes must be comparable (in a sense that will be specified later). When applying our methodology to real datasets, we found that the sample sizes were far from balanced, with some being very large. Another point is that the computation of the variance of the test statistic increases substantially with the sample sizes. To deal with these two issues, we adopt a strategy which parallels to that in \cite{Cuesta2010}. Specifically, balanced random samples are taken from the data and the associated $p$-value is obtained. This is repeated several times and the $p$-values are adjusted.




The remainder of the paper is organized as follows. Section~\ref{section:teststatistic} introduces the proposed test statistic for comparing general parameters across multiple populations. Section~\ref{section:asymptoticdist} develops the asymptotic distribution of the test statistic under the null hypothesis. Section~\ref{section:variance} presents a ratio-consistent estimator for the asymptotic variance. Section~\ref{section:power} studies the asymptotic power of the proposed test under various types of alternatives. 
Section~\ref{section:computation} discusses computational challenges and proposed solutions.

\section{Test statistic}
\label{section:teststatistic}

Let \( X_{1}, \dots, X_{k} \) be independent random elements defined on a common probability space \( (\Omega, \mathcal{F}, \Prob) \), taking values in a measurable space \( (E, \mathcal{E}) \), where $E$ is an arbitrary space and $\mathcal{E}$ is $\sigma$-algebra on $E$. 

Furthermore, let \( m \geq 1 \) be a fixed integer. We consider that  
\begin{equation}
\begin{aligned}
    &\bm{X}_i = \{X_{i1}, \dots, X_{in_i} \} \text{ is a random sample of size } n_i \geq 2m \text{ from the distribution of } X_i, \\
    &\text{and the samples } \bm{X}_1, \ldots, \bm{X}_k \text{ are independent.}
\end{aligned}
\label{eq:independent_sample}
\end{equation}
Let \( \theta_i \in \mathbb{R} \) be defined as  
\[
\theta_i = \E\left \{h(X_{i1},\ldots, X_{im})\right\}, \; \text{for}\; i=1,\ldots k,
\]  
where \( h: E^m \to \mathbb{R} \) is a symmetric kernel of degree $m$.  For this kernel, we denote its variance terms as  
\begin{equation}  
    \sigma^2_{ij} = \V\{h_j(X_{i1},\ldots,X_{ij})\}, \quad j=1,\ldots,m,
    \label{eq:definitionsigmai}  
\end{equation}  
where $h_j(x_1,\ldots,x_j)=\E\{h(X_{i1},\ldots, X_{im})|X_{i1}=x_1,\ldots,X_{ij}=x_j\}$.

Our objective is to test $H_0$ when $k$ is large.
Under the null, this is equivalent to testing whether \( D_k = 0 \), where  
\[
D_k = \frac{1}{k} \sum_{i=1}^k (\theta_i - \overline{\theta})^2 = \left( \frac{1}{k}-\frac{1}{k^2}\right) \sum_{i=1}^k \theta_i^2  -  \frac{1}{k^2}  \sum_{i \neq j} \theta_{i} \theta_{j},
\]
with \( \overline{\theta} = \frac{1}{k} \sum_{i=1}^k \theta_i \).  

Given the definition of the parameter \( \theta_i \),  
it is natural to estimate it using the \( U \)-statistic  
\[
\widehat{\theta}_i = \frac{1}{n_i(m)} \sum_{j_1 \neq \cdots \neq j_m} h(X_{ij_1}, \ldots, X_{ij_m}), 
\]
where \( n_i(k) = n_i(n_i-1) \cdots (n_i-k+1) \), which  
is an unbiased estimator of $\theta_i$ for $i=1,\ldots,k$. Similarly, to obtain an unbiased estimator of \( D_k \), we first construct a unbiased estimators of $ \theta_i^2$  for $i=1,\ldots,k$ by applying \( U \)-statistics theory. The kernel associated with \( \theta^2 \) is  
\begin{equation}
    H(x_{1},\ldots, x_{2m})= \frac{1}{(2m)!} \sum_{j_1\neq \cdots \neq j_{2m}} h(x_{j_1},\ldots,x_{j_m}) h(x_{j_{m+1}},\ldots,x_{j_{2m}}),
    \label{eq:H}
\end{equation}  
and the corresponding unbiased estimator of \( \theta_i^2 \) is given by  
\[
\widehat{\theta_i^2} = \frac{1}{n_i(2m)} \sum_{j_1\neq \cdots \neq j_{2m}} h(X_{ij_1},\ldots,X_{ij_m}) h(X_{ij_{m+1}},\ldots,X_{ij_{2m}}), \; \text{for}\; i=1,\ldots,k.
\]
Then, we define the statistic  
\begin{equation}
T_k = \left( \frac{1}{k} - \frac{1}{k^2} \right) \sum_{i=1}^k \widehat{\theta_i^2} - \frac{1}{k^2} \sum_{i \neq j} \widehat{\theta}_i \widehat{\theta}_j,
\label{eq:Tk}
\end{equation}
which is, by construction, an unbiased estimator of \( D_k \). Moreover, under certain conditions, \( T_k \) is a consistent estimator of \( D_k \), meaning that \( T_k \) converges in probability to \( D_k \). Specifically, we assume that  
\begin{equation}
    \E\{h^2(X_{i1}, \ldots, X_{im})\} < M, \quad \text{for some positive constant } M \text{ and } \forall i,
    \label{ass:boundenessh2}
\end{equation}
which is a key assumption in most of our results to derive meaningful conclusions. Throughout this paper, all limits are understood to be taken as $k \to \infty$, unless otherwise specified.

\begin{lemma}
\label{lemma:consistency}
Suppose that \eqref{eq:independent_sample} and \eqref{ass:boundenessh2} hold. Then,  
\[
T_k - D_k \tp 0 .
\]
Moreover, if \(\E\{|h(X_{i1},\ldots,X_{im})|^{2+\delta}\} < M\), $\forall i $  for some \(M > 0\), then
\[
T_k- D_k \tas 0.
\]
\end{lemma}
The null hypothesis is identified by \( D_k = 0 \), whereas the alternative hypothesis is characterized by \( D_k > 0 \). Since \( T_k \) is unbiased and consistent for \( D_k \), it would be reasonable to reject \( H_0 \) when \( T_k \) is large. However, to determine when \( T_k \) is large, we need to know its distribution under \( H_0 \). Since we have made no assumptions about the underlying distribution of \( X_1, \ldots, X_k \), obtaining the exact distribution of \( T_k \) is not feasible. For this reason, we approximate it through its asymptotic distribution, which is derived in the next section. 

\section{Asymptotic distribution }
\label{section:asymptoticdist}

To derive the asymptotic distribution of \( T_k \), we decompose the statistic as stated in the following lemma:
\begin{lemma}
\label{lemma:Tkdecomp}
$T_k$ can be expressed as
\begin{equation}
T_k = D_k + T_{k,\text{Lin}} + R_k,
\label{eq:decompTk}
\end{equation}
where
\[
T_{k,\text{Lin}} = \sum_{i=1}^{k} \frac{L_{ik}}{k}, \quad \text{ with }\quad L_{ik} = \widehat{\theta_i^2} - \theta^2_i - 2\overline{\theta}(\widehat{\theta}_i - \theta_i)
\]
and
\[
R_k = - \frac{1}{k^2} \sum_{i=1}^{k} \left( \widehat{\theta_i^2} - \theta^2_i \right) + \frac{2}{k^2} \sum_{i=1}^{k} \theta_i (\widehat{\theta}_i - \theta_i) - \frac{1}{k^2} \sum_{i \neq j} (\widehat{\theta}_i - \theta_i)(\widehat{\theta}_j - \theta_j).
\]
\end{lemma}
Note that the variable $L_{ik}$ depend only on the sample $\bm{X}_i$. As the samples $\bm{X}_1, \ldots, \bm{X}_k$ are independent, the variables $L_{1k}, \ldots, L_{kk}$ are also independent. Consequently,
\begin{equation}
    \V(T_{k,\text{Lin}}) = 
     \frac{1}{k} V_k, \qquad V_K=\frac{1}{k} \sum_{i=1}^{k} \V(L_{ik}).
    \label{eq:vtklin}
\end{equation}

\begin{lemma}
\label{lemma:variancegeneral}
Assume that $\E\{h^2(X_{i1},\ldots,X_{im})\}< \infty$. The variance of $L_{ik}$ is given by
\begin{equation*}
    \V(L_{ik})=\V(\mathbb{L}_{i1})+\Big(\theta_i-\overline{\theta}\Big)^2\V(2\widehat{\theta}_i)+\Big(\theta_i-\overline{\theta}\Big)\Xi^3_{i},
\end{equation*}
where $\mathbb{L}_{i1}= \widehat{\theta^2} - \theta_i^2 - 2\theta_i(\widehat{\theta_i} - \theta_i).$
\end{lemma}
Now, we outline the scenarios in which we will obtain the asymptotic distribution of $T_k$. We consider two cases. The first case is the obvious one, where $H_0$ holds. If $H_0$ is not true, we assume one of the following
assumptions describing the behavior of the sequence \(\{D_k\}\). 

\begin{assumption}
\(H_0\) does not hold, one of these conditions is satisfied:
\begin{enumerate}[label=(\roman*)]
    \item \(D_k \geq \delta, \forall k > k_0\), for some \(\delta > 0\) and some \(k_0 \in \mathbb{N}\).
    \item $D_k \to 0$ and $a_k D_k \geq \delta, \forall k > k_0$ for some $\delta>0$ and some $k_0 \in \mathbb{N}$, where $\{a_k\}$ is an increasing sequence such that $a_k \to \infty$, $0< a_k\leq k^{\eta}$, for some $0<\eta<1$, and $a_k(\theta_i-\overline{\theta})^2 \leq M$, $\forall i$, for some positive constant $M$.
\end{enumerate}
\label{assumption}
\end{assumption}
Condition (i) corresponds to a scenario where \(D_k\) remains bounded away from zero for all sufficiently large \(k\). In contrast, condition (ii) allows \(D_k\) to approach zero at a controlled rate $a_k$, since
\begin{equation}
    \delta\leq a_k D_k= \frac{1}{k} \sum_{i=1}^k a_k(\theta_i-\overline{\theta})^2\leq M.
    \label{eq:controlDk}
\end{equation} 
These two conditions are assumed to ensure that the sequence $\{D_k\}$ behaves regularly, allowing us to derive meaningful outcomes.

To establish the asymptotic distribution of the test statistic, we will show that the variance of the remainder term is negligible compared to the variance of the linear term. To achieve this, we make two additional assumptions.

First, we assume that the sample sizes are comparable, in the sense that
\begin{equation}
    n_i(k) = c_i(k) n_0(k), \quad n_0(k) \geq 2m, \quad 1 \leq c_i(k) \leq C, \text{ for some positive constant } C
    \label{eq:compsamp}
\end{equation}
The comparability condition on sample sizes, as defined in \eqref{eq:compsamp}, is imposed to mitigate the impact of severe sample size imbalances, which affects the asymptotic distribution of the statistic. Such assumption is commonly made in the statistical literature. For instance, \cite{scholz1987k} and \cite{RizzoSzekely2010} impose similar conditions to establish the consistency of the test developed in their respective works. In what follows, we omit the explicit dependence on  $k$  of  $n_0$  and  $c_i$  for simplicity in notation.

Next, we assume that $\widehat{\theta_i}$ is non-degenerate, meaning that $\sigma_{i1}^2 > 0$, where $\sigma_{i1}^2$ is defined in \eqref{eq:definitionsigmai}. 
Morover, to derive the next result we assume a uniform lower bound on $\sigma_{i1}^2$, i.e.,
\begin{equation}
\sigma^2_{i1} \geq \rho , \; \text{for some positive constant }\rho \text{ and } \forall i.
    \label{ass:boundsigma}
\end{equation}

\begin{lemma}
    \label{lemma:VRk0}
Assume that \eqref{eq:independent_sample},\eqref{ass:boundenessh2}, \eqref{eq:compsamp} and \eqref{ass:boundsigma} hold. Also, assume that either \(H_0\) holds or Assumption \ref{assumption} is satisfied. Then, \(\V(R_k)/\V(T_{k,\text{Lin}}) \to 0\).
\end{lemma}

Once we have the previous result, deriving the asymptotic distribution of $T_k$ reduces to determining the asymptotic distribution of $T_{k,\text{Lin}}$, which is stated in the following result.
\begin{lemma}
Assume that the conditions of Lemma \ref{lemma:VRk0} are satisfied, and that if Assumption \ref{assumption}.ii holds, then either \( a_k/n_0 \) is bounded or \( a_k/n_0 \to \infty \). Under these conditions, we have that \( T_{k, \text{Lin}} / \sqrt{\V(T_{k, \text{Lin}})} \tl Z \), where \( Z \) follows a standard normal distribution.
\label{lemma:Tlinlimit}
\end{lemma}
\begin{theorem}
Under the assumptions of Lemma \ref{lemma:Tlinlimit}, we have that $(T_k - D_k) / \sqrt{\V(T_{k, \text{Lin}})} \tl Z$. 
\label{Theorem:Tklim2}
\end{theorem}

Theorem \ref{Theorem:Tklim2} holds either if \( H_0 \) is true or Assumption \ref{assumption} is satisfied. Specifically, when \( H_0 \) holds, we have \( D_k = 0 \), and from \eqref{eq:vtklin}, it follows that \( \sqrt{k} T_k / \sqrt{V_k} \tl Z \). As discussed at the end of Section \ref{section:teststatistic}, it is reasonable to reject the null hypothesis whenever \( T_k \) is large, as larger values of \( T_k \) indicate a stronger deviation from \( H_0 \). Theorem \ref{Theorem:Tklim2} provides a rigorous foundation for constructing statistical tests. Specifically, for a chosen significance level \( \alpha \in (0, 1) \), we can reject \( H_0 \) if \( \sqrt{k} T_k / \sqrt{V_k} > z_{1-\alpha} \), where \( z_{1-\alpha} \) is the critical value corresponding to the \( (1-\alpha) \)-quantile of the standard normal distribution. This decision rule ensures that the probability of incorrectly rejecting \( H_0 \) (the type I error) does not exceed \( \alpha \) for large $k$.

However, in practice, Theorem \ref{Theorem:Tklim2} cannot be directly applied because the variance \( V_k \) is typically unknown, as no assumptions are made about the distribution of the random variables. To address this limitation, we require a ratio-consistent estimator \( \widehat{V}_k \) of \( V_k \) under \( H_0 \). Specifically, this means that \( \widehat{V}_k / V_k \tp 1 \) when \( H_0 \) holds. By Slutsky's theorem, substituting \( V_k \) with \( \widehat{V}_k \) in the normalization does not alter the limiting distribution. Consequently, constructing such an estimator is a crucial step for the practical implementation of the test. The details of this construction are presented in the following section.

\begin{remark}
    Note that all previous results are valid whenever the sample sizes satisfy \eqref{eq:compsamp}, which allows the sample size to remain bounded or increase with $k$.
\end{remark}
\section{Variance estimation}
\label{section:variance}

The variance we are trying to estimate, as defined in \eqref{eq:vtklin}, is given by
\[
V_k = \frac{1}{k} \sum_{i=1}^k \V(L_{ik}).
\]
Firstly, we eliminate the aditive constants in $\V(L_{ik})$, as they have no effect in the variance computation. For that, we introduce the auxiliar quantities 
\begin{equation}
    l_{ik} = \widehat{\theta_i^2} - 2\overline{\theta} \widehat{\theta_i},
    \label{eq:lik}
\end{equation} and by construction, we have that \( \E[l_{ik}] = \theta_i^2 - 2\overline{\theta} \theta_i \)  for $i=1,\ldots,k$. Observe that \( L_{ik} = l_{ik} - \E[l_{ik}] \), and therefore, \( \V(L_{ik}) = \V(l_{ik}) \). Let
\[S_{l}^2=\frac{1}{k-1} \sum_{i=1}^k(l_{ik}-\overline{l})^2, \quad \overline{l}=\frac{1}{k}\sum_{i=1}^k l_{ik}.\]
Its expetation is given by
\begin{align*}
    \E(S_l^2)&=\frac{1}{k-1}\E\left\{ \sum_{i=1}^k (l_{ik}^2-2\overline{l}l_{ik}+\overline{l}^2)\right\} \\ 
    &=\frac{1}{k-1}\Big[ k\{V_k + \overline{\E^2(l)}\}-2\{V_k+k\overline{\E(l)}^2\}+\{V_k+k\overline{\E[l]}^2\}\Big] \\
    &=V_k+ \frac{k}{k-1}\left\{\overline{\E^2(l)}-\overline{\E(l)}^2\right\},
\end{align*}
where $\overline{\E^2(l)}=\frac{1}{k}\sum_{i=1}^k \E^2(l_{ik})$ and $\overline{\E(l)}^2=\Big\{\frac{1}{k}\sum_{i=1}^k \E(l_{ik})\Big\}^2$. We have that, \[\overline{\E^2(l)}-\overline{\E(l)}^2=\frac{1}{k}\sum_{i=1}^k\left\{\theta_i^2-\overline{\theta^2}-2\overline{\theta}(\theta_i-\overline{\theta})\right\}^2,
\]
with $\overline{\theta^2}=\frac{1}{k}\sum_{i=1}^k \theta_i^2$. Therefore we can write
\[\E(S_l^2)=V_k+\Delta_k,\]
where
\begin{equation}
    \Delta_k=\frac{1}{k-1}\sum_{i=1}^k\left\{\theta_i^2-\overline{\theta^2}-2\overline{\theta}(\theta_i-\overline{\theta})\right\}^2.
    \label{eq:Deltak}
\end{equation}

It is easy to observe that, under \( H_0 \), \(\Delta_k = 0\), and hence \(\E(S_l^2) = V_k\). However, although \( S_l^2 \) could serve as a ratio-consistent estimator of \(\E(S_l^2)\), it cannot be directly computed because \( l_{1k}, \ldots, l_{kk} \) are unknown, as \(\overline{\theta}\) is unknown. To address this issue, we introduce a slight modification. Define
\[
\hat{l}_{ik} = \widehat{\theta_i^2} - 2 \widehat{\overline{\theta}}\widehat{\theta_i}, \quad i = 1, \ldots, k, \quad \widehat{\overline{\theta}} = \frac{1}{k} \sum_{\substack{j=1}}^k \widehat{\theta}_j,
\]
and 
\begin{equation}
  S_{\hat{l}}^2 = \frac{1}{k-1} \sum_{i=1}^k (\hat{l}_{ik} - \overline{\hat{l}})^2, \quad \overline{\hat{l}} = \frac{1}{k} \sum_{i=1}^k \hat{l}_{ik}.
  \label{eq:Shatl}
\end{equation}

Finally, we show that \( S^2_{\hat{l}} \) is a ratio-consistent estimator of \(\E(S^2_l)\), as formally stated in the following result.

\begin{proposition}
Suppose that the assumptions in Lemma \ref{lemma:Tlinlimit} hold, that $n_0/k \to 0$, 
and that for some $\delta>0$, one of the following cases holds:
\begin{enumerate}
    \item $H_0$ is true and $\E\left(\left|n_0 L_{ik}\right|^{2+\delta}\right) \leq M, \, \forall i,$ 
    \item Assumption \ref{assumption}.i is true and $\E\left(\left|\sqrt{n_0} L_{ik}\right|^{2+\delta}\right) \leq M, \, \forall i,$ 
    \item Assumption \ref{assumption}.ii is true, $a_k/n_0 \leq M$ and $\E\left(\left|\sqrt{a_k n_0} L_{ik}\right|^{2+\delta}\right) \leq M, \, \forall i,$ 
    \item Assumption \ref{assumption}.ii is true, $a_k/n_0 \to \infty$ and $\E\left(\left|n_0 L_{ik}\right|^{2+\delta}\right) \leq M, \, \forall i.$ 
\end{enumerate}
Then $S^2_{\hat{l}}/\E(S^2_l) \tp 1$.
    \label{proposition:varianceestimator}
\end{proposition}

In contrast to the results previously stated, the ratio-consistency of \( S^2_{\hat{l}} \) in Proposition~\ref{proposition:varianceestimator} additionally requires that \( n_0/k \to 0 \). This assumption is not particularly restrictive in practice, as it is generally unrealistic to obtain increasing amounts of data from all populations as the number of populations grows. Furthermore, the proposition assumes the existence of a $(2+\delta)$-th moment for some $\delta > 0$. This condition is justified in each of the cases considered, since Lemma~\ref{lemma:boundsVgen} (parts 1, 3, 4, and 5) in Section~\ref{section: proofs} shows that the corresponding variables, $n_0 L_{ik}$, $\sqrt{n_0} L_{ik}$, or $\sqrt{a_k n_0} L_{ik}$, have finite second moments. In light of this, it is reasonable to impose a slightly stronger moment condition in each setting.

Observe that under \(H_0\), since \(\Delta_k = 0\), the variance estimator satisfies the ratio-consistency property, i.e., $S^2_{\hat{l}}/\E(S^2_l) \tp 1$. Recalling Theorem \ref{Theorem:Tklim2} and applying Slutsky’s theorem, it follows that  
\begin{equation}
    \mathcal{T}_k=\sqrt{k}\frac{ T_k}{S_{\hat{l}}} \tl Z.
    \label{eq:statistic}
\end{equation}
Consequently, the test that rejects \( H_0 \) if  $\mathcal{T}_k > z_{1-\alpha},$ has level $\alpha$ asymptotically.

\section{Asymptotic power}
\label{section:power}

In this section, we study the asymptotic power of the proposed test under the alternatives outlined in Assumption~\ref{assumption}. The following theorem summarizes the cases in which the test is consistent, 
i.e., $\mathcal{P}\to 1$, where $\mathcal{P}=\lim \Prob(\mathcal{T}_k > z_{1-\alpha})$.

\begin{theorem}
Assume that the conditions in Proposition~\ref{proposition:varianceestimator} are satisfied, and that one of the following scenarios holds:
\begin{enumerate}[label=\alph*), leftmargin=*]  
\item Assumption~\ref{assumption}.i is satisfied.  
\item Assumption~\ref{assumption}.ii is satisfied and \( \sqrt{k}/a_k \to \infty \).  
\item Assumption~\ref{assumption}.ii is satisfied, \( M_1\leq \sqrt{k}/a_k\leq M_2\), and \( n_0 \to \infty \).  
\item Assumption~\ref{assumption}.ii is satisfied, \( \sqrt{k}/a_k \to 0 \), \( n_0 \to \infty \), and \( a_k/n_0 \) remains bounded.  
\item Assumption~\ref{assumption}.ii is satisfied, \( \sqrt{k}/a_k \to 0 \), \( n_0 \to \infty \), \( a_k/n_0 \to \infty \), and \( \sqrt{k}n_0/a_k \to \infty \).
\end{enumerate}  
Then, 
$\mathcal{P}=1$.

\label{theorem:power}
\end{theorem}

The following results examine the power of the test under alternative scenarios, 
where the asymptotic power is less than $1$.

\begin{proposition}
\begin{enumerate}
    \item Assume the conditions of case (c) in Theorem \ref{theorem:power}, but with \( n_0 \) bounded instead of tending to infinity. Then, $\alpha < \mathcal{P} < 1$. Furthermore, if \( \sqrt{k} D_k / \sqrt{V_k} \to d \in (0, \infty) \), then
    $\mathcal{P} = \Phi(z_{\alpha} + d).$

    
    \item Assume the conditions of case (e) in Theorem \ref{theorem:power}, but with \(M_1\leq \sqrt{k}n_0/a_k \leq M_2 \). Then,
    $\alpha < \mathcal{P} < 1.
    $ Furthermore, if \( \sqrt{k} D_k / \sqrt{V_k} \to d \in (0, \infty) \), then
      $\mathcal{P}=\Phi(z_\alpha+d).$

    \item Assume the conditions of case (e) in Theorem \ref{theorem:power}, but with \( \sqrt{k}n_0/a_k \to 0 \). Then, $\mathcal{P} = \alpha.$
    
\end{enumerate}
\label{proposition:conditioncfixedn0}
\end{proposition}

The results presented above are based on the possible behaviors of the sequences \( \{n_0\} \) and \( \{a_k\} \). These include cases where \( n_0 \) is bounded or tends to infinity, as well as different rates at which \( n_0 \) and \( a_k \) approach infinity. They provide intuitive insights into the asymptotic power of the test, which can be summarized as follows:
\begin{itemize}
\item Alternatives satisfying Assumption~\ref{assumption}.i, where \( D_k \) is bounded away from zero, can be detected by the test both when \( n_0 \) is bounded and when \( n_0 \to \infty \) (case~(a)).

\item For alternatives satisfying Assumption~\ref{assumption}.ii, the behavior depends on the rate \( a_k \) at which \( D_k \) converges to zero. The following subcases can be distinguished:
\begin{itemize}
\item If \( a_k = o(\sqrt{k}) \) (case~(b)), \( D_k \) converges to 0 at a slow rate, allowing the test to detect the violation of \( H_0 \) both when \( n_0 \) is bounded and when \( n_0 \to \infty \).

\item If \( a_k \sim \sqrt{k} \) (i.e., \( M_1 \sqrt{k} \leq a_k \leq M_2 \sqrt{k}\); case~(c) and Proposition \ref{proposition:conditioncfixedn0}.1 , \( D_k \) converges faster than in the previous case, and a bounded \( n_0 \) is no longer sufficient. Consistency is achieved when \( n_0 \to \infty \). On the other hand, if \( n_0 \) remains bounded, the power converges to a value strictly between \( \alpha \) and 1. This value depends on the limiting behavior of the signal-to-noise ratio \( \sqrt{k} D_k / \sqrt{V_k} \), and is equal to \( \Phi(z_\alpha + d) \) when this quantity converges to \( d \in (0, \infty) \).

\item For the final cases, where \( \sqrt{k} = o(a_k) \) (cases~(d) and~(e) in Theorem \ref{theorem:power}, and parts~2--3 of Proposition \ref{proposition:conditioncfixedn0}, \( D_k \) converges to zero so fast that the consistency of the test is achieved if \( n_0 \to \infty \) at a rate faster than a certain threshold. In particular, consistency is achieved when \( n_0 \geq M a_k \). Otherwise, the test is consistent when \( n_0 \) grows faster than \( a_k / \sqrt{k} \). Finally, when \( n_0 \sim a_k / \sqrt{k} \), the power tends to a value in \( (\alpha, 1) \), which again depends on the limiting signal-to-noise ratio. If instead \( n_0 = o(a_k / \sqrt{k}) \), then the power tends to \( \alpha \).
\end{itemize}
\end{itemize}

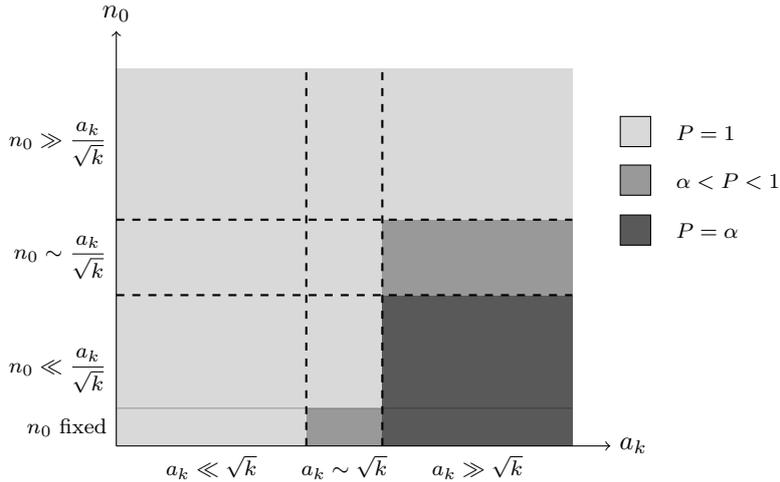
\begin{figure}
\centering
\begin{tikzpicture}[scale=1]

  \def\xa{0}
  \def\xb{2.5}
  \def\xc{3.5}
  \def\xd{6.0}
  \def\ya{0}
  \def\yb{0.5}
  \def\yc{2.0}
  \def\yd{3.0}
  \def\ye{5.0}
  \def\xf{6.5}
  \def\yf{5.5}

  \definecolor{colP1}{gray}{0.85}    
  \definecolor{colPinter}{gray}{0.6} 
  \definecolor{colPalpha}{gray}{0.35}

  \fill[colP1]     (\xa,\ya) rectangle (\xb,\yb); 
  \fill[colP1]     (\xa,\yb) rectangle (\xb,\yc); 
  \fill[colP1]     (\xa,\yc) rectangle (\xb,\yd); 
  \fill[colP1]     (\xa,\yd) rectangle (\xb,\ye); 

  \fill[colPinter] (\xb,\ya) rectangle (\xc,\yb); 
  \fill[colP1]     (\xb,\yb) rectangle (\xc,\yc); 
  \fill[colP1]     (\xb,\yc) rectangle (\xc,\yd); 
  \fill[colP1]     (\xb,\yd) rectangle (\xc,\ye); 

  \fill[colPalpha] (\xc,\ya) rectangle (\xd,\yb); 
  \fill[colPalpha] (\xc,\yb) rectangle (\xd,\yc); 
  \fill[colPinter] (\xc,\yc) rectangle (\xd,\yd); 
  \fill[colP1]     (\xc,\yd) rectangle (\xd,\ye); 

  \draw[->] (0,0) -- (\xf,0) node[right] {$a_k$};
  \draw[->] (0,0) -- (0,\yf) node[above] {$n_0$};

  \draw[opacity=0.3, thin] (0,\yb) -- (\xd,\yb); 
  \draw[dashed, thick] (0,\yc) -- (\xd,\yc); 
  \draw[dashed, thick] (0,\yd) -- (\xd,\yd); 
  \draw[dashed, thick] (\xb,0) -- (\xb,\ye);
  \draw[dashed, thick] (\xc,0) -- (\xc,\ye);

  \node[left] at (0,1) {\footnotesize $n_0 \ll \dfrac{a_k}{\sqrt{k}}$};
  \node[left] at (0,0.25) {\footnotesize $n_0$ fixed};
  \node[left] at (0,2.5) {\footnotesize $n_0 \sim \dfrac{a_k}{\sqrt{k}}$};
  \node[left] at (0,4) {\footnotesize $n_0 \gg \dfrac{a_k}{\sqrt{k}}$};

  \node[below] at (1.25,0) {\footnotesize $a_k \ll \sqrt{k}$};
  \node[below] at (3.0,0)   {\footnotesize $a_k \sim \sqrt{k}$};
  \node[below] at (4.75,0)  {\footnotesize $a_k \gg \sqrt{k}$};

  \matrix[anchor=north west, row sep=6pt, column sep=6pt] at (\xd+0.5,\yf-1.0) {
    \node[draw, fill=colP1,    minimum size=4mm, inner sep=0pt] {}; & 
    \node[align=left] {\footnotesize $P = 1$}; \\
    
    \node[draw, fill=colPinter, minimum size=4mm, inner sep=0pt] {}; & 
    \node[align=left] {\footnotesize $\alpha < P < 1$}; \\
    
    \node[draw, fill=colPalpha, minimum size=4mm, inner sep=0pt] {}; & 
    \node[align=left] {\footnotesize $P = \alpha$}; \\
  };

\end{tikzpicture}

\caption{Scheme of the asymptotic power of the test for different rates of $a_k$ and $n_0$.}
\label{graph}
\end{figure}

All above cases are graphically displayed in Figure \ref{graph}.

\noindent

This discussion shows that the more difficult it is to detect the violation of \( H_0 \), the larger the number of observations in each group must be to ensure consistency. While easier alternatives can be detected with a bounded \( n_0 \), more challenging ones require not only \( n_0 \to \infty \), but also at a fast enough rate.

\section{Computational considerations}
\label{section:computation}



The practical application of the proposed methodology may face two main challenges. First, the behavior of the test statistic \( \mathcal{T}_k \) may deviate from its asymptotic distribution, particularly when the number of groups \( k \) is small or moderate. Second, when group sizes are large or highly unbalanced, both the computational burden and the validity of the theoretical approximations may be adversely affected. To address these challenges, we propose two complementary strategies:
\begin{itemize}
    \item A linear bootstrap procedure to improve the practical performance of the test for small to moderate $k$.
    \item A random sampling approach that reduces the computational burden and mitigates the impact of group size imbalance.
\end{itemize}

\subsection{Linear bootstrap}

In the simulation study, the convergence to the normal law of the null distribution of the test statistic \( \mathcal{T}_k \) was found to be relatively slow for small to moderate $k$. To improve the performance of the test, we employed a linear bootstrap procedure based on resampling the linear component of the statistic, as proposed in \cite{JimenezGamero2025}. The procedure is implemented as follows:
\begin{enumerate}
    \item Compute the centered linear components \( \hat{l}_{1k} - \bar{\hat{l}}, \ldots, \hat{l}_{kk} - \bar{\hat{l}} \).
    \item Generate a bootstrap sample \( L^*_{1k}, \ldots, L^*_{kk} \) from the empirical distribution of the centered values.
    \item Evaluate the linear approximation of \( \mathcal{T}_k \) on the bootstrap sample, obtaining \( \mathcal{T}_k^* = \sqrt{k}\,\overline{L^*} / \sqrt{\widehat{V}^*} \), where \( \overline{L^*} = (1/k) \sum_{i=1}^k L^*_{ik} \) and \( \widehat{V}^* = \{1/(k-1)\} \sum_{i=1}^k (L^*_{ik} - \overline{L^*})^2 \).
    \item Approximate the null distribution of \( \mathcal{T}_k \) by the conditional distribution of \( \mathcal{T}_k^* \) given the observed data \( \bm{X}_1,\ldots,\bm{X}_k \).
\end{enumerate}

In practice, this approximation is computed via Monte Carlo simulation: steps 2 and 3 are repeated \( B \) times, yielding bootstrap replicates \( \mathcal{T}_k^{*1}, \ldots, \mathcal{T}_k^{*B} \). The empirical distribution of these replicates provides an estimate of the linear bootstrap distribution of \( \mathcal{T}_k \).

This resampling strategy is particularly efficient for large \( n_1,\ldots,n_k \), as it completely avoids the computation of \( \widehat{\theta^2_1}, \ldots, \widehat{\theta^2_k} \) at bootstrap level, which are the most computationally demanding components of the statistic.

\subsection{Random sampling}
\label{sect:6.2}
Real datasets often exhibit considerable imbalance in sample sizes, with some groups having disproportionately large sizes. This fact poses three issues: a) to derive the results, comparability of sample sizes is assumed (see \eqref{eq:compsamp}); b) the sample sizes must be small relative to $k$ (recall that Proposition \ref{proposition:varianceestimator} assumes that $n_0/k \to 0$); c) the estimation of the variance requires computing \( \widehat{\theta^2_i} \), each based on a \( U \)-statistic of degree \( 2m \) with computational cost \( \mathcal{O}(n_i^{2m}) \) for \( i = 1, \ldots, k \), and hence, the overall computational burden may become substantially large. 

To address these limitations, we propose a random sampling scheme whereby, for each application of the test, we draw a sample of size \( \tilde{n} \ll k \) from each group. This reduces the computational burden and provides a better approximation to the condition \( n_0 / k \to 0 \) in Proposition \ref{proposition:varianceestimator}. It is worth noting that different random samples may lead to contradictory conclusions. To mitigate this problem, we adopt the strategy proposed by~\cite{Cuesta2010} in the context of random projection tests for functional data. Their approach consists in generating \( L \) independent random projections, computing the \( p \)-value for each of them, and applying a multiple testing correction (such as the method proposed in~\cite{BenjaminiYekutieli2001}). 
Specifically, their methodology computes \( L \) \( p \)-values, one for each random projection; then sorts them as \( p_{(1)} \leq \cdots \leq p_{(L)} \); and finally obtains the  \( p \)-value as $p_{\text{adj}} = \min_{1 \leq i \leq L} \left\{  (L/i) p_{(i) } \right\}$. The global null hypothesis is then rejected at level \( \alpha \) whenever \( p_{\text{adj}} \leq \alpha \). Proposition 2.3 in \cite{Cuesta2010} ensures that the level is at most $\alpha$. 

In our setting, we apply the same idea by drawing \( L \) independent random samples from each group; calculating the $p$-value for each set of random samples and calculating $p_{\text{adj}}$ as before. The authors recommend the use up to $L = 30$.
Furthermore, as shown in Theorem~\ref{theorem:power}, some alternatives require \( n_0 \) to increase with \( k \) in order to ensure the consistency of the test. Consequently, rather than fixing a small sample size in advance, we adopt an adaptive strategy in which \( \tilde{n} \) is progressively increased until the test decision stabilizes.


\section{Simulation studies} 
\label{section:simulation}

In this section, we evaluate the finite-sample performance of the proposed test through Monte Carlo simulations selecting a particular choice of the kernel $h$. We have selected two different statistics to evaluate the performance of the test: the Gini mean difference (GMD) and the Spearman's rho $(\rho_S)$. We assess the performance of two tests:
\begin{itemize}
    \item The test based on the asymptotic normal approximation derived in this work (N).
    \item The linear bootstrap test discussed in Section \ref{section:computation} (LB).
\end{itemize}

To reduce the computational burden involved in bootstrap, we employ the warp-speed method proposed in \cite{Giacomini2013}. In this approach, instead of computing critical values for each Monte Carlo sample, we generate a single resample per iteration. Let \( J \) denote the total number of Monte Carlo replications. In each iteration \( j = 1, \dots, J \), the test statistic is first computed on the original sample, resulting in \( \mathcal{T}_{k,j} \), and then on a resample generated in the same iteration, yielding \( \mathcal{T}_{k,j}^* \). The empirical distribution of the resampled statistics, \( \{ \mathcal{T}_{k,1}^*, \dots, \mathcal{T}_{k,J}^* \} \), is then used to compute $p$-values for the corresponding original values \( \mathcal{T}_{k,j} \), $j=1,\ldots,J$.

For each configuration, we generate independent samples from \( k \) populations, with identical sample sizes \( n_i = n_0 \). The $p$-value of the observed value of the test statistic is computed using both the normal approximation and the bootstrap approach. The simulation is repeated $J=10,000$ times to estimate the empirical rejection rates under the null and alternative hypotheses. We study the impact of both the number of populations \( k \) and the individual sample size \( n_0 \) on the performance of the test.

\subsection{Gini mean difference}
\label{subsection:gmd}

The Gini mean difference (GMD) of a random variable \( X \) is defined as \( \text{GMD} = \mathbb{E}(|X_1 - X_2|) \), where \( X_1 \) and \( X_2 \) are independent copies of \( X \). This corresponds to a $U$-statistic of degree 2 with kernel \( h(x_1, x_2) = |x_1 - x_2| \).

\subsubsection*{Type I error: Experimental design and results}

To empirically asses the probability of type I error, we considered four data-generating settings. In the first three, all populations are identically distributed, following:
\begin{enumerate}
    \item a standard normal distribution \( \mathcal{N}(0,1) \),
    \item a chi-squared distribution with 3 degrees of freedom, \( \chi^2(3) \),
    \item a Student’s \( t \)-distribution with 5 degrees of freedom \( t(5) \),
    \item A heterogeneous mixture: one third of the populations are sampled from each of the aforementioned distributions, and the data are rescaled so that all populations have a Gini mean difference equal to 1.
\end{enumerate}
 For each case, we considered \( k \in \{50, 100, 500, 1000, 1500\} \) populations and equal group sizes \( n_i \in \{5, 10, 15\}\) for $i=1,\ldots,k$. The results are reported in Table \ref{table:levelgini}.

\begin{table}[h]
\centering
\footnotesize
\renewcommand{\arraystretch}{1.1}
\newcommand{\pair}[2]{\makebox[1.5em][r]{#1}\hspace{1em}\makebox[1.5em][r]{#2}}
\begin{tabular}{@{}cc@{}}
  \begin{minipage}[t]{0.44\textwidth}
    \centering
    {\bfseries $\mathcal{N}(0,1)$}\\[1ex]
    \begin{tabular}{l c c c c}
      \toprule
      $k$ &  & $n_i=5$ & $n_i=10$ & $n_i=15$ \\
      \midrule
      \multirow{2}{*}{50}
        & N & \pair{4.4}{9.5} & \pair{3.2}{7.7} & \pair{3.1}{7.4} \\
        & LB & \pair{4.6}{9.4} & \pair{4.7}{9.6} & \pair{4.7}{9.9} \\
      \midrule
      \multirow{2}{*}{100}
        & N & \pair{4.5}{9.5} & \pair{3.4}{8.4} & \pair{3.4}{8.4} \\
        & LB & \pair{5.0}{10.1}& \pair{4.3}{9.2} & \pair{5.0}{10.2} \\
      \midrule
      \multirow{2}{*}{500}
        & N & \pair{4.5}{9.4} & \pair{4.3}{8.8} & \pair{4.0}{9.1} \\
        & LB & \pair{4.5}{9.4} & \pair{5.3}{10.0}& \pair{5.0}{10.3} \\
      \midrule
      \multirow{2}{*}{1000}
        & N & \pair{4.6}{9.9} & \pair{4.3}{9.3} & \pair{4.2}{9.0} \\
        & LB & \pair{4.8}{9.9} & \pair{5.0}{10.1}& \pair{5.0}{10.1} \\
      \midrule
      \multirow{2}{*}{1500}
        & N & \pair{4.7}{9.8} & \pair{4.3}{8.9} & \pair{4.3}{9.1} \\
        & LB & \pair{4.9}{10.1}& \pair{4.6}{9.8} & \pair{4.8}{9.6} \\
      \bottomrule
    \end{tabular}
  \end{minipage}%
  \hspace{0.02\textwidth}%
  \begin{minipage}[t]{0.44\textwidth}
    \centering
    {\bfseries $\chi^2(3)$}\\[1ex]
    \begin{tabular}{l c c c c}
      \toprule
      $k$ &  & $n_i=5$ & $n_i=10$ & $n_i=15$ \\
      \midrule
      \multirow{2}{*}{50}
        & N & \pair{4.9}{9.9} & \pair{4.0}{9.0} & \pair{3.6}{8.8} \\
        & LB & \pair{3.6}{8.7} & \pair{3.9}{9.3} & \pair{4.1}{10.0} \\
      \midrule
      \multirow{2}{*}{100}
        & N & \pair{4.8}{10.1}& \pair{4.0}{9.1} & \pair{3.8}{8.2} \\
        & LB & \pair{4.3}{9.7} & \pair{4.1}{9.4} & \pair{4.3}{8.5} \\
      \midrule
      \multirow{2}{*}{500}
        & N & \pair{4.9}{9.9} & \pair{4.6}{9.4} & \pair{4.2}{8.8} \\
        & LB & \pair{4.7}{9.6} & \pair{4.8}{9.8} & \pair{4.7}{9.9} \\
      \midrule
      \multirow{2}{*}{1000}
        & N & \pair{4.4}{9.6} & \pair{4.6}{9.6} & \pair{4.6}{9.6} \\
        & LB & \pair{5.0}{9.8} & \pair{4.8}{9.6} & \pair{4.7}{10.3} \\
      \midrule
      \multirow{2}{*}{1500}
        & N & \pair{4.7}{9.9} & \pair{4.8}{9.7} & \pair{4.7}{9.9} \\
        & LB & \pair{4.7}{10.0}& \pair{5.1}{9.8} & \pair{5.0}{10.1} \\
      \bottomrule
    \end{tabular}
  \end{minipage}

  \\[48ex] 

  \begin{minipage}[t]{0.44\textwidth}
    \centering
    {\bfseries $t(5)$}\\[1ex]
    \begin{tabular}{l c c c c}
      \toprule
      $k$ & Test & $n_i=5$ & $n_i=10$ & $n_i=15$ \\
      \midrule
      \multirow{2}{*}{50}
        & N & \pair{3.9}{8.8} & \pair{3.2}{8.1} & \pair{3.4}{8.1} \\
        & LB & \pair{3.3}{8.3} & \pair{3.6}{8.2} & \pair{4.5}{9.3} \\
      \midrule
      \multirow{2}{*}{100}
        & N & \pair{4.1}{8.8} & \pair{3.2}{8.2} & \pair{3.4}{8.2} \\
        & LB & \pair{4.0}{8.8} & \pair{3.8}{8.9} & \pair{4.5}{10.2} \\
      \midrule
      \multirow{2}{*}{500}
        & N & \pair{4.1}{9.0} & \pair{3.8}{9.3} & \pair{3.9}{8.5} \\
        & LB & \pair{5.0}{9.6} & \pair{4.5}{9.9} & \pair{4.8}{9.9} \\
      \midrule
      \multirow{2}{*}{1000}
        & N & \pair{4.4}{9.4} & \pair{4.2}{9.2} & \pair{4.2}{9.2} \\
        & LB & \pair{4.7}{9.8} & \pair{5.0}{9.9} & \pair{5.1}{9.6} \\
      \midrule
      \multirow{2}{*}{1500}
        & N & \pair{4.6}{9.9} & \pair{4.2}{9.2} & \pair{4.3}{9.4} \\
        & LB & \pair{4.9}{10.3}& \pair{5.2}{9.8} & \pair{5.0}{10.3} \\
      \bottomrule
    \end{tabular}
  \end{minipage}%
  \hspace{0.02\textwidth}%
  \begin{minipage}[t]{0.44\textwidth}
    \centering
    { \phantom{t(5)} $ \text{Mixture}$ \phantom{t(5)}}\\[1ex]
    \begin{tabular}{l c c c c}
      \toprule
      $k$ &  & $n_i=5$ & $n_i=10$ & $n_i=15$ \\
      \midrule
      \multirow{2}{*}{50}
        & N & \pair{4.5}{9.8} & \pair{3.7}{8.7} & \pair{3.2}{8.1} \\
        & LB & \pair{3.5}{8.8} & \pair{3.7}{9.0} & \pair{3.9}{8.6} \\
      \midrule
      \multirow{2}{*}{100}
        & N & \pair{4.3}{9.3} & \pair{4.0}{8.9} & \pair{3.5}{8.3} \\
        & LB & \pair{4.1}{9.4} & \pair{4.4}{9.8} & \pair{4.2}{9.6} \\
      \midrule
      \multirow{2}{*}{500}
        & N & \pair{4.4}{9.4} & \pair{4.3}{9.4} & \pair{4.3}{9.2} \\
        & LB & \pair{4.4}{9.2} & \pair{5.2}{10.1}& \pair{5.1}{10.4} \\
      \midrule
      \multirow{2}{*}{1000}
        & N & \pair{4.3}{9.2} & \pair{4.4}{9.4} & \pair{4.2}{9.1} \\
        & LB & \pair{4.8}{9.4} & \pair{5.2}{9.9} & \pair{4.8}{9.6} \\
      \midrule
      \multirow{2}{*}{1500}
        & N & \pair{4.7}{9.7} & \pair{4.1}{9.4} & \pair{4.3}{9.3} \\
        & LB & \pair{5.2}{10.3}& \pair{4.7}{10.2}& \pair{4.6}{9.7} \\
      \bottomrule
    \end{tabular}
  \end{minipage}
\end{tabular}

\caption{Empirical Type I error rates expressed in percentages, for nominal levels $\alpha=5\%$ (left) and $\alpha=10\%$ (right), for the normal approximation (N) and the linear bootstrap approximation (LB).}
\label{table:levelgini}
\end{table}
Looking at Table~\ref{table:levelgini}, we observe that that both the asymptotic test (N) and the linear bootstrap test (LB) tend to be conservative for small $k$. For small values of $k$, increasing the sample size $n_i$ can lead to a decrease in the empirical type I error rate, specially for the asymptotic test. As the number of populations increases, both procedures converge to the nominal levels. In general, the  bootstrap test display a better performance.

\subsubsection*{Empirical power: Experimental design and results}

We generated scaled normal data with $\text{GMD}=1$ and created alternatives by assigning a proportion $\pi$ of contaminated data with  \( \pi \in \{0.1, 0.2, 0.3, 0.4\} \) of the populations with GMD equal to \( \theta \in \{1.25, 1.5\} \). We considered configurations with \( k \in \{50, 100, 200\} \) and \( n_0 \in \{5, 10, 15\} \). In this case, we have that \( D_k = \pi(1 - \pi)(\theta-1)^2 \). Therefore, by varying \( \pi \) and \( \theta \), we assess the power of the test under alternatives located at different distances from the null. The results obtained are reported in Table~\ref{tab:powergini}.

\begin{table}[h]
\centering
\renewcommand{\arraystretch}{1}
\scriptsize
\newcommand{\pair}[2]{%
  \makebox[1.5em][r]{#1}%
  \hspace{1em}%
  \makebox[1.5em][r]{#2}%
}
\resizebox{\textwidth}{!}{%
  \begin{tabular}{@{}llc ccc ccc@{}}
    \toprule
    \multirow{2}{*}{$\pi$} & \multirow{2}{*}{$k$} & \multirow{2}{*}{Test}
      & \multicolumn{3}{c}{$\theta=1.25$}
      & \multicolumn{3}{c}{$\theta=1.5$} \\
    \cmidrule(lr){4-6} \cmidrule(lr){7-9}
      & & & $n_i=5$ & $n_i=10$ & $n_i=15$
      & $n_i=5$ & $n_i=10$ & $n_i=15$ \\
    \midrule

    \multirow{6}{*}{0.1}
      & 50  & N  
        & \pair{5.5}{11.0} & \pair{5.9}{12.7} & \pair{7.8}{16.5}
        & \pair{7.5}{14.6} & \pair{13.4}{25.6} & \pair{23.5}{40.5} \\
      &     & LB 
        & \pair{5.7}{10.9} & \pair{8.7}{15.7} & \pair{11.2}{20.7}
        & \pair{7.7}{14.9} & \pair{18.4}{31.0} & \pair{30.9}{47.0} \\
    \cmidrule(lr){2-9}
      & 100 & N  
        & \pair{5.8}{11.9} & \pair{8.1}{16.6} & \pair{11.8}{22.8}
        & \pair{9.3}{17.4} & \pair{24.0}{38.7} & \pair{44.5}{62.3} \\
      &     & LB 
        & \pair{6.4}{12.8} & \pair{10.5}{18.9} & \pair{16.1}{28.0}
        & \pair{10.4}{18.7} & \pair{28.1}{42.6} & \pair{53.6}{68.9} \\
    \cmidrule(lr){2-9}
      & 200 & N  
        & \pair{6.8}{13.3} & \pair{11.1}{21.2} & \pair{19.0}{32.2}
        & \pair{12.7}{22.5} & \pair{41.1}{57.6} & \pair{74.0}{85.7} \\
      &     & LB 
        & \pair{7.6}{14.7} & \pair{13.6}{23.8} & \pair{22.9}{35.9}
        & \pair{14.2}{24.1} & \pair{45.3}{60.8} & \pair{78.5}{87.8} \\
    \midrule

    \multirow{6}{*}{0.2}
      & 50  & N  
        & \pair{6.4}{12.6} & \pair{9.2}{18.2}  & \pair{14.9}{27.6}
        & \pair{11.5}{20.9} & \pair{33.0}{50.2} & \pair{61.6}{77.9} \\
      &     & LB 
        & \pair{6.4}{12.8} & \pair{12.6}{22.6} & \pair{20.8}{33.0}
        & \pair{11.6}{21.3} & \pair{40.4}{56.0} & \pair{69.0}{82.6} \\
    \cmidrule(lr){2-9}
      & 100 & N  
        & \pair{7.3}{14.2} & \pair{14.6}{25.9} & \pair{25.7}{41.1}
        & \pair{16.0}{27.1} & \pair{57.9}{73.1} & \pair{89.8}{95.6} \\
      &     & LB 
        & \pair{8.2}{15.1} & \pair{17.2}{28.4} & \pair{32.5}{47.1}
        & \pair{17.8}{28.6} & \pair{62.2}{76.6} & \pair{93.1}{96.5} \\
    \cmidrule(lr){2-9}
      & 200 & N  
        & \pair{9.2}{17.1} & \pair{22.7}{36.9} & \pair{44.5}{60.9}
        & \pair{25.3}{38.8} & \pair{85.5}{92.3} & \pair{99.6}{99.9} \\
      &     & LB 
        & \pair{10.5}{18.4} & \pair{26.4}{40.4} & \pair{51.0}{64.7}
        & \pair{27.6}{40.7} & \pair{87.5}{93.4} & \pair{99.8}{99.9} \\
    \midrule

    \multirow{6}{*}{0.3}
      & 50  & N  
        & \pair{7.1}{14.3} & \pair{12.3}{23.5} & \pair{22.4}{37.7}
        & \pair{15.7}{26.4} & \pair{53.1}{69.7} & \pair{86.1}{93.7} \\
      &     & LB 
        & \pair{7.6}{14.7} & \pair{16.6}{28.3} & \pair{29.3}{44.2}
        & \pair{16.4}{27.5} & \pair{60.3}{74.2} & \pair{89.5}{95.5} \\
    \cmidrule(lr){2-9}
      & 100 & N  
        & \pair{8.7}{16.4} & \pair{21.1}{34.5} & \pair{40.1}{56.0}
        & \pair{23.1}{36.7} & \pair{81.2}{90.2} & \pair{98.9}{99.7} \\
      &     & LB 
        & \pair{9.5}{17.6} & \pair{24.6}{37.9} & \pair{47.5}{61.0}
        & \pair{25.5}{38.9} & \pair{84.4}{91.7} & \pair{99.4}{99.8} \\
    \cmidrule(lr){2-9}
      & 200 & N  
        & \pair{11.6}{20.9} & \pair{35.1}{50.8} & \pair{66.1}{79.1}
        & \pair{38.8}{54.0} & \pair{97.8}{99.3} & \pair{100}{100} \\
      &     & LB 
        & \pair{13.1}{22.4} & \pair{40.9}{54.7} & \pair{71.9}{81.7}
        & \pair{42.6}{56.0} & \pair{98.5}{99.5} & \pair{100}{100} \\
    \midrule

    \multirow{6}{*}{0.4}
      & 50  & N  
        & \pair{7.5}{14.9} & \pair{15.1}{27.4} & \pair{28.6}{45.3}
        & \pair{18.6}{30.7} & \pair{65.8}{79.9} & \pair{94.3}{97.9} \\
      &     & LB 
        & \pair{8.0}{15.3} & \pair{19.8}{32.5} & \pair{36.5}{51.9}
        & \pair{20.7}{31.7} & \pair{72.7}{83.3} & \pair{96.3}{98.5} \\
    \cmidrule(lr){2-9}
      & 100 & N  
        & \pair{9.5}{17.9} & \pair{26.1}{40.6} & \pair{49.7}{65.7}
        & \pair{29.5}{44.2} & \pair{91.4}{96.1} & \pair{99.9}{100} \\
      &     & LB 
        & \pair{10.7}{19.1} & \pair{30.2}{44.9} & \pair{56.9}{70.3}
        & \pair{33.2}{47.0} & \pair{93.2}{96.8} & \pair{100}{100} \\
    \cmidrule(lr){2-9}
      & 200 & N  
        & \pair{13.4}{23.2} & \pair{43.9}{60.2} & \pair{78.1}{88.1}
        & \pair{49.6}{64.6} & \pair{99.7}{99.9} & \pair{100}{100} \\
      &     & LB 
        & \pair{14.8}{24.1} & \pair{49.4}{64.1} & \pair{82.9}{89.9}
        & \pair{53.3}{66.0} & \pair{99.8}{100} & \pair{100}{100} \\
    \bottomrule
  \end{tabular}%
}
\caption{Empirical power in percentages for normal distributed data, with \(100\pi\%\) of the observations having \(\text{GMD}=\theta\) and \(100(1-\pi)\%\) having $\text{GMD}=1$, with nominal levels  $\alpha=5\%$ (left) and $\alpha=10\%$ (right), for the normal approximation (N) and the linear bootstrap approximation (LB).}
\label{tab:powergini}
\end{table}

As can be seen in Table~\ref{tab:powergini}, the empirical power increases with higher values of \( \pi \), \( \theta \), and \( n_0 \). The bootstrap approximation slightly outperforms the asymptotic approach.

\subsection{Spearman's rho}
\label{subsection:spearman}

The Spearman’s rho between the \( r \)-th and \( s \)-th components of a \( d \)-dimensional random vector \( X = (X^{(1)}, \ldots, X^{(d)}) \) is defined as
\[
\rho_{S}^{(r,s)} = \mathbb{E}[h^{(r,s)}(X_{1}, X_{2}, X_{3})],
\]
where \( X_1, X_2 \), and \( X_3 \) are independent copies of \( X \), and \( h^{(r,s)} \) is a symmetric kernel of degree 3 given by
\[
h^{(r,s)}(x_1, x_2, x_3) = \frac{1}{2} \sum_{1 \leq \alpha \neq \beta \neq \gamma \leq 3} \text{sgn}(x_\alpha^{(r)} - x_\beta^{(r)}) \cdot \text{sgn}(x_\alpha^{(s)} - x_\gamma^{(s)}).
\]
The superscripts \((r)\) and \((s)\) indicate the \( r \)-th and \( s \)-th components of the vectors \( x_1, x_2, x_3 \in \mathbb{R}^d \), respectively. This \( U \)-statistic representation can be found in \cite{ElMaache2003}

\subsubsection*{Type I error: Experimental design and results}

To empirically assess the probability of type I error of the test when the parameter of interest is Spearman’s rho, we generated data from bivariate normal distributions with varying degrees of correlation between two components. Specifically, we considered Spearman’s rho values \( \rho_S \in \{0, 0.25, 0.5, 0.75\} \),  \( k \in \{50, 100, 500, 1000, 1500, 2000\} \) populations, and group sizes \( n_0 \in \{7, 10, 15\} \).
\begin{table}[h]
\centering
\footnotesize
\renewcommand{\arraystretch}{1.1}
\newcommand{\pair}[2]{\makebox[1.5em][r]{#1}\hspace{1em}\makebox[1.5em][r]{#2}}
\begin{tabular}{@{}cc@{}}

  \begin{minipage}[t]{0.48\textwidth}
    \centering
    {\bfseries $\rho_S=0$}\\[1ex]
    \begin{tabular}{l c c c c}
      \toprule
      $k$ & Test & $n_i=7$ & $n_i=10$ & $n_i=15$ \\
      \midrule
      \multirow{2}{*}{50}
        & N & \pair{4.3}{9.4} & \pair{3.3}{7.8} & \pair{3.0}{7.6} \\
        & LB & \pair{4.6}{10.2}& \pair{5.1}{9.7} & \pair{4.8}{9.3} \\
      \midrule
      \multirow{2}{*}{100}
        & N & \pair{4.3}{9.4} & \pair{4.0}{8.9} & \pair{3.4}{8.0} \\
        & LB & \pair{5.1}{9.9}& \pair{4.9}{10}& \pair{4.8}{9.8} \\
      \midrule
      \multirow{2}{*}{500}
        & N & \pair{4.8}{9.7} & \pair{4.4}{9.5} & \pair{3.9}{8.4} \\
        & LB & \pair{4.9}{10.4}& \pair{4.8}{9.8} & \pair{4.8}{9.3}\\
      \midrule
      \multirow{2}{*}{1000}
        & N & \pair{4.7}{9.8} & \pair{4.5}{9.6} & \pair{4.3}{9.4} \\
        & LB & \pair{4.7}{10.1}& \pair{5.2}{10.2} & \pair{4.9}{10.1}\\
      \midrule
      \multirow{2}{*}{1500}
        & N & \pair{4.5}{9.6} & \pair{4.7}{9.6} & \pair{4.6}{9.5} \\
        & LB & \pair{4.8}{9.8} & \pair{5.3}{10.2}& \pair{5.4}{10.2}\\
      \midrule
      \multirow{2}{*}{2000}
        & N & \pair{4.5}{9.5} & \pair{4.7}{9.8} & \pair{4.7}{9.7} \\
        & LB & \pair{4.7}{10.0} & \pair{5.2}{10.4}& \pair{5.3}{10.5}\\
      \bottomrule
    \end{tabular}
  \end{minipage}%
  \hspace{0.02\textwidth}%
  \begin{minipage}[t]{0.48\textwidth}
    \centering
    {\bfseries $\rho_S=0.25$}\\[1ex]
    \begin{tabular}{l c c c c}
      \toprule
      $k$ & Test & $n_i=7$ & $n_i=10$ & $n_i=15$ \\
      \midrule
      \multirow{2}{*}{50}
        & N & \pair{3.8}{8.7} & \pair{2.9}{7.3} & \pair{2.8}{7.2} \\
        & LB & \pair{4.4}{9.7} & \pair{4.2}{8.9} & \pair{4.7}{9.5} \\
      \midrule
      \multirow{2}{*}{100}
        & N & \pair{3.9}{8.3} & \pair{3.2}{7.8} & \pair{3.0}{7.8} \\
        & LB & \pair{4.7}{9.8} & \pair{4.8}{9.6} & \pair{4.5}{9.6} \\
      \midrule
      \multirow{2}{*}{500}
        & N & \pair{4.2}{9.2} & \pair{3.8}{8.5} & \pair{3.9}{9.0} \\
        & LB & \pair{4.6}{9.2} & \pair{4.5}{9.8} & \pair{4.7}{10.4}\\
      \midrule
      \multirow{2}{*}{1000}
        & N & \pair{4.5}{9.6} & \pair{4.4}{8.9} & \pair{4.2}{8.8} \\
        & LB & \pair{4.7}{9.5} & \pair{5.1}{9.6} & \pair{4.9}{9.4} \\
      \midrule
      \multirow{2}{*}{1500}
        & N & \pair{4.6}{9.6} & \pair{4.3}{9.3} & \pair{4.0}{8.7} \\
        & LB & \pair{4.7}{10.0}& \pair{4.7}{10.2}& \pair{4.6}{10.0}\\
      \midrule
      \multirow{2}{*}{2000}
        & N & \pair{4.6}{9.6} & \pair{4.6}{9.2} & \pair{4.0}{8.9} \\
        & LB & \pair{4.9}{9.9} & \pair{5.2}{9.9} & \pair{4.6}{9.8} \\
      \bottomrule
    \end{tabular}
  \end{minipage}

  \\[55ex] 

  \begin{minipage}[t]{0.48\textwidth}
    \centering
    {\bfseries $\rho_S=0.5$}\\[1ex]
    \begin{tabular}{l c c c c}
      \toprule
      $k$ & Test & $n_i=7$ & $n_i=10$ & $n_i=15$ \\
      \midrule
      \multirow{2}{*}{50}
        & N & \pair{2.5}{6.7} & \pair{2.4}{6.5} & \pair{2.4}{6.5} \\
        & LB & \pair{3.4}{8.6} & \pair{3.8}{8.4} & \pair{3.8}{8.7} \\
      \midrule
      \multirow{2}{*}{100}
        & N & \pair{2.8}{7.2} & \pair{2.9}{7.3} & \pair{2.6}{6.9} \\
        & LB & \pair{4.2}{9.0} & \pair{5.0}{9.8} & \pair{4.0}{8.9} \\
      \midrule
      \multirow{2}{*}{500}
        & N & \pair{3.6}{8.3} & \pair{3.4}{8.2} & \pair{3.5}{8.5} \\
        & LB & \pair{4.5}{9.3} & \pair{4.5}{9.5} & \pair{4.4}{9.8}\\
      \midrule
      \multirow{2}{*}{1000}
        & N & \pair{4.2}{8.7} & \pair{4.0}{0.1} & \pair{4.0}{8.6} \\
        & LB & \pair{4.8}{9.6} & \pair{5.2}{10.2}& \pair{4.7}{9.4} \\
      \midrule
      \multirow{2}{*}{1500}
        & N & \pair{4.2}{9.1} & \pair{4.3}{9.2} & \pair{4.2}{8.6} \\
        & LB & \pair{4.7}{9.5} & \pair{5.2}{10.3}& \pair{5.0}{9.9} \\
      \midrule
      \multirow{2}{*}{2000}
        & N & \pair{4.4}{9.6} & \pair{4.5}{9.2} & \pair{4.3}{8.9} \\
        & LB & \pair{4.5}{10.0}& \pair{5.3}{10.2}& \pair{5.4}{10.1}\\
      \bottomrule
    \end{tabular}
  \end{minipage}%
  \hspace{0.02\textwidth}%
  \begin{minipage}[t]{0.48\textwidth}
    \centering
    {\bfseries $\rho_S=0.75$}\\[1ex]
    \begin{tabular}{l c c c c}
      \toprule
      $k$ & Test & $n_i=7$ & $n_i=10$ & $n_i=15$ \\
      \midrule
      \multirow{2}{*}{50}
        & N & \pair{1.1}{4.3} & \pair{1.7}{4.9} & \pair{1.6}{5.4} \\
        & LB & \pair{1.6}{5.9} & \pair{2.5}{6.5} & \pair{3.4}{7.7} \\
      \midrule
      \multirow{2}{*}{100}
        & N & \pair{1.4}{5.1} & \pair{1.8}{5.2} & \pair{1.9}{5.3} \\
        & LB & \pair{2.9}{7.6} & \pair{3.4}{8.2} & \pair{3.1}{7.1} \\
      \midrule
      \multirow{2}{*}{500}
        & N & \pair{2.3}{6.9} & \pair{2.7}{7.0} & \pair{2.9}{7.3} \\
        & LB & \pair{4.6}{9.5} & \pair{4.0}{8.7} & \pair{4.4}{9.4} \\
      \midrule
      \multirow{2}{*}{1000}
        & N & \pair{2.9}{7.9} & \pair{3.0}{7.4} & \pair{3.2}{8.0} \\
        & LB & \pair{4.4}{10.0}& \pair{4.3}{9.5} & \pair{4.5}{9.6} \\
      \midrule
      \multirow{2}{*}{1500}
        & N & \pair{3.2}{7.9} & \pair{3.4}{8.0} & \pair{3.8}{8.0} \\
        & LB & \pair{4.8}{9.9} & \pair{4.7}{10.3}& \pair{5.2}{9.8} \\
      \midrule
      \multirow{2}{*}{2000}
        & N & \pair{3.4}{7.9} & \pair{3.5}{8.1} & \pair{3.5}{8.4} \\
        & LB & \pair{4.8}{9.4} & \pair{5.3}{10.1}& \pair{5.0}{10.1}\\
      \bottomrule
    \end{tabular}
  \end{minipage}

\end{tabular}
\caption{Empirical Type I error rates expressed in percentages, for bivariate normal data with varying $\rho_S$ values with nominal levels $\alpha=5\%$ (left) and $\alpha=10\%$ (right), for the normal approximation (N) and the linear bootstrap approximation (LB).}
\label{table:levelspearman}
\end{table}
The results are summarized in Table~\ref{table:levelspearman}. It can be observed that the empirical levels of the bootstrap test are generally closer to the nominal values than those of the asymptotic test. The discrepancy between both tests diminishes as $k$ increases. For small values of $k$ and $\rho_S > 0$, both procedures tend to be conservative, particularly as $\rho_S$ becomes larger.


\subsubsection*{Power: Experimental design and results}

To empirically study the power, we generated standard bivariate normal data and created alternatives by assigning a proportion \( \pi \in \{0.1, 0.2, 0.3, 0.4\} \) of the populations a nonzero Spearman’s rho equal to \( \theta \in \{0.25, 0.5\} \), while the remaining populations kept \( \rho_S = 0 \). We used \( k \in \{50, 100, 200\} \) and \( n_0 \in \{7, 10, 15\} \). In this case, we have that \( D_k = \pi(1 - \pi)\theta^2 \), which equals zero under the null hypothesis. By varying \( \pi \) and \( \theta \), we assess the power of the test under alternatives located at different distances from the null. The results obtaind are presented in Table~\ref{tab:powerrhoS}. As in the GMD case, the power increases with $\pi, \theta, n_0$ and $k$.

\begin{table}[ht]
\centering
\renewcommand{\arraystretch}{1}
\scriptsize
\newcommand{\pair}[2]{%
  \makebox[1.5em][r]{#1}%
  \hspace{1em}%
  \makebox[1.5em][r]{#2}%
}
\resizebox{\textwidth}{!}{%
  \begin{tabular}{@{}llc ccc ccc@{}}
    \toprule
    \multirow{2}{*}{$\pi$} & \multirow{2}{*}{$k$} & \multirow{2}{*}{Test}
      & \multicolumn{3}{c}{$\theta=0.25$}
      & \multicolumn{3}{c}{$\theta=0.5$} \\
    \cmidrule(lr){4-6} \cmidrule(lr){7-9}
      & & & $n_i=7$ & $n_i=10$ & $n_i=15$
      & $n_i=7$ & $n_i=10$ & $n_i=15$ \\
    \midrule

    \multirow{6}{*}{0.1} & \multirow{2}{*}{50} & N  
      & \pair{5.3}{10.7} & \pair{4.8}{11.1} & \pair{5.5}{12.8}
      & \pair{8.5}{17.3} & \pair{13.4}{27.4} & \pair{25.5}{46.4} \\
      &                      & LB 
      & \pair{6.4}{12.1} & \pair{6.9}{13.3} & \pair{8.6}{16.8}
      & \pair{11.3}{20.6} & \pair{22.0}{34.2} & \pair{38.6}{56.7} \\
    \cmidrule(lr){2-9}
    & \multirow{2}{*}{100} & N  
      & \pair{5.5}{11.5} & \pair{6.1}{13.1} & \pair{8.5}{16.9}
      & \pair{13.0}{24.4} & \pair{25.9}{42.7} & \pair{52.8}{71.9} \\
      &                      & LB 
      & \pair{6.7}{12.9} & \pair{7.9}{14.8} & \pair{11.3}{19.9}
      & \pair{17.5}{28.1} & \pair{34.5}{49.1} & \pair{65.3}{78.5} \\
    \cmidrule(lr){2-9}
    & \multirow{2}{*}{200} & N  
      & \pair{6.5}{13.2} & \pair{8.2}{16.1} & \pair{13.1}{23.5}
      & \pair{21.4}{35.3} & \pair{47.7}{64.9} & \pair{85.7}{93.7} \\
      &                      & LB 
      & \pair{7.5}{13.9} & \pair{10.0}{17.9} & \pair{16.0}{26.1}
      & \pair{25.7}{39.5} & \pair{55.2}{69.2} & \pair{90.3}{95.4} \\
    \midrule

    \multirow{6}{*}{0.2} & \multirow{2}{*}{50} & N  
      & \pair{6.2}{12.3} & \pair{6.6}{14.1} & \pair{9.3}{19.0}
      & \pair{16.0}{28.5} & \pair{32.1}{49.9} & \pair{64.2}{80.8} \\
      &                      & LB 
      & \pair{7.5}{13.7} & \pair{9.0}{16.7} & \pair{13.5}{23.4}
      & \pair{19.8}{32.0} & \pair{41.1}{56.8} & \pair{75.1}{85.7} \\
    \cmidrule(lr){2-9}
    & \multirow{2}{*}{100} & N  
      & \pair{7.2}{14.3} & \pair{14.6}{25.9} & \pair{25.7}{41.1}
      & \pair{16.0}{27.1} & \pair{57.9}{73.1} & \pair{89.8}{95.6} \\
      &                      & LB 
      & \pair{8.2}{15.1} & \pair{17.2}{28.4} & \pair{32.5}{47.1}
      & \pair{17.8}{28.6} & \pair{62.2}{76.6} & \pair{93.1}{96.5} \\
    \cmidrule(lr){2-9}
    & \multirow{2}{*}{200} & N  
      & \pair{9.2}{17.1} & \pair{22.7}{36.9} & \pair{44.5}{60.9}
      & \pair{25.3}{38.8} & \pair{85.5}{92.3} & \pair{99.6}{99.9} \\
      &                      & LB 
      & \pair{10.5}{18.4} & \pair{26.4}{40.4} & \pair{51.0}{64.7}
      & \pair{27.6}{40.7} & \pair{87.5}{93.4} & \pair{99.8}{99.9} \\
    \midrule

    \multirow{6}{*}{0.3} & \multirow{2}{*}{50} & N  
      & \pair{6.6}{13.2} & \pair{12.3}{23.5} & \pair{22.4}{37.7}
      & \pair{15.7}{26.4} & \pair{53.1}{69.7} & \pair{86.1}{93.7} \\
      &                      & LB 
      & \pair{7.6}{14.7} & \pair{16.6}{28.3} & \pair{29.3}{44.2}
      & \pair{16.4}{27.5} & \pair{60.3}{74.2} & \pair{89.5}{95.5} \\
    \cmidrule(lr){2-9}
    & \multirow{2}{*}{100} & N  
      & \pair{8.7}{16.4} & \pair{21.1}{34.5} & \pair{40.1}{56.0}
      & \pair{23.1}{36.7} & \pair{81.2}{90.2} & \pair{98.9}{99.7} \\
      &                      & LB 
      & \pair{9.5}{17.6} & \pair{24.6}{37.9} & \pair{47.5}{61.0}
      & \pair{25.5}{38.9} & \pair{84.4}{91.7} & \pair{99.4}{99.8} \\
    \cmidrule(lr){2-9}
    & \multirow{2}{*}{200} & N  
      & \pair{11.6}{20.9} & \pair{35.1}{50.8} & \pair{66.1}{79.1}
      & \pair{38.8}{54.0} & \pair{97.8}{99.3} & \pair{100}{100} \\
      &                      & LB 
      & \pair{13.1}{22.4} & \pair{40.9}{54.7} & \pair{71.9}{81.7}
      & \pair{42.6}{56.0} & \pair{98.5}{99.5} & \pair{100}{100} \\
    \midrule

    \multirow{6}{*}{0.4} & \multirow{2}{*}{50} & N  
      & \pair{7.1}{14.9} & \pair{15.1}{27.4} & \pair{28.6}{45.3}
      & \pair{18.6}{30.7} & \pair{65.8}{79.9} & \pair{94.3}{97.9} \\
      &                      & LB 
      & \pair{8.0}{15.3} & \pair{19.8}{32.5} & \pair{36.5}{51.9}
      & \pair{20.7}{31.7} & \pair{72.7}{83.3} & \pair{96.3}{98.5} \\
    \cmidrule(lr){2-9}
    & \multirow{2}{*}{100} & N  
      & \pair{9.5}{17.9} & \pair{26.1}{40.6} & \pair{49.7}{65.7}
      & \pair{29.5}{44.2} & \pair{91.4}{96.1} & \pair{99.9}{100} \\
      &                      & LB 
      & \pair{10.7}{19.1} & \pair{30.2}{44.9} & \pair{56.9}{70.3}
      & \pair{33.2}{47.0} & \pair{93.2}{96.8} & \pair{100}{100} \\
    \cmidrule(lr){2-9}
    & \multirow{2}{*}{200} & N  
      & \pair{13.4}{23.2} & \pair{43.9}{60.2} & \pair{78.1}{88.1}
      & \pair{49.6}{64.6} & \pair{99.7}{99.9} & \pair{100}{100} \\
      &                      & LB 
      & \pair{14.8}{24.1} & \pair{49.4}{64.1} & \pair{82.9}{89.9}
      & \pair{53.3}{66.0} & \pair{99.8}{100} & \pair{100}{100} \\
    \bottomrule
  \end{tabular}%
}
\caption{Empirical power in percentages for standard normally distributed data, with \(100\pi\%\) of the observations having $\rho_S=\theta$ and \(100(1-\pi)\%\) having $\rho_S=0$. For  $\alpha=5\%$ (left) and $\alpha=10\%$ (right), for the test based on the normal approximation (N) and the linear bootstrap approximation (LB).}
\label{tab:powerrhoS}
\end{table}

\section{Empirical application}
\label{section:application}

In this section, we illustrate the application of the proposed test using data from the IPUMS USA Version~16.0 dataset \cite{ruggles2025ipums}. IPUMS USA Version~16.0 is a harmonized microdata extract that combines decennial census Public Use Microdata Samples (PUMS) and American Community Survey (ACS) 1-year samples from 2001 to 2023. This dataset provides consistent individual-level information on demographics (e.g., age, sex, race), education, employment, income, housing, and health, along with geographic identifiers down to the Public Use Microdata Area (PUMA) and county levels for sub-state analysis.

For our analysis, we select the ACS 2023 1-year sample from IPUMS USA Version~16.0. This specific sample allows us to conduct a detailed cross-sectional examination of income and employment patterns across counties. To construct a meaningful partition of the data, we combine two key geographic identifiers: \texttt{PUMA}, which represents Public Use Microdata Areas, and \texttt{COUNTYICP}, a coding scheme for counties based on the Inter-University Consortium for Political and Social Research (ICPSR) standard. The concatenation of these two variables uniquely identifies each subregion in the dataset. The study focuses on two main variables: \texttt{INCEARN}, which corresponds to the individual's earned income, and \texttt{UHRSWORK}, representing the usual number of hours worked per week. After filtering out records with missing, unknown, or non-identifiable values for the relevant variables, the final sample consists of 1,132,892 individuals distributed across $k=1{,}587$ counties. The sample sizes at the county level range from 193 to 2,691 individuals



\subsection*{Gini mean difference of income relative to weekly hours worked}\label{subsec:total_income}

\begin{figure}[ht]
    \centering
    \includegraphics[width=0.7\linewidth]{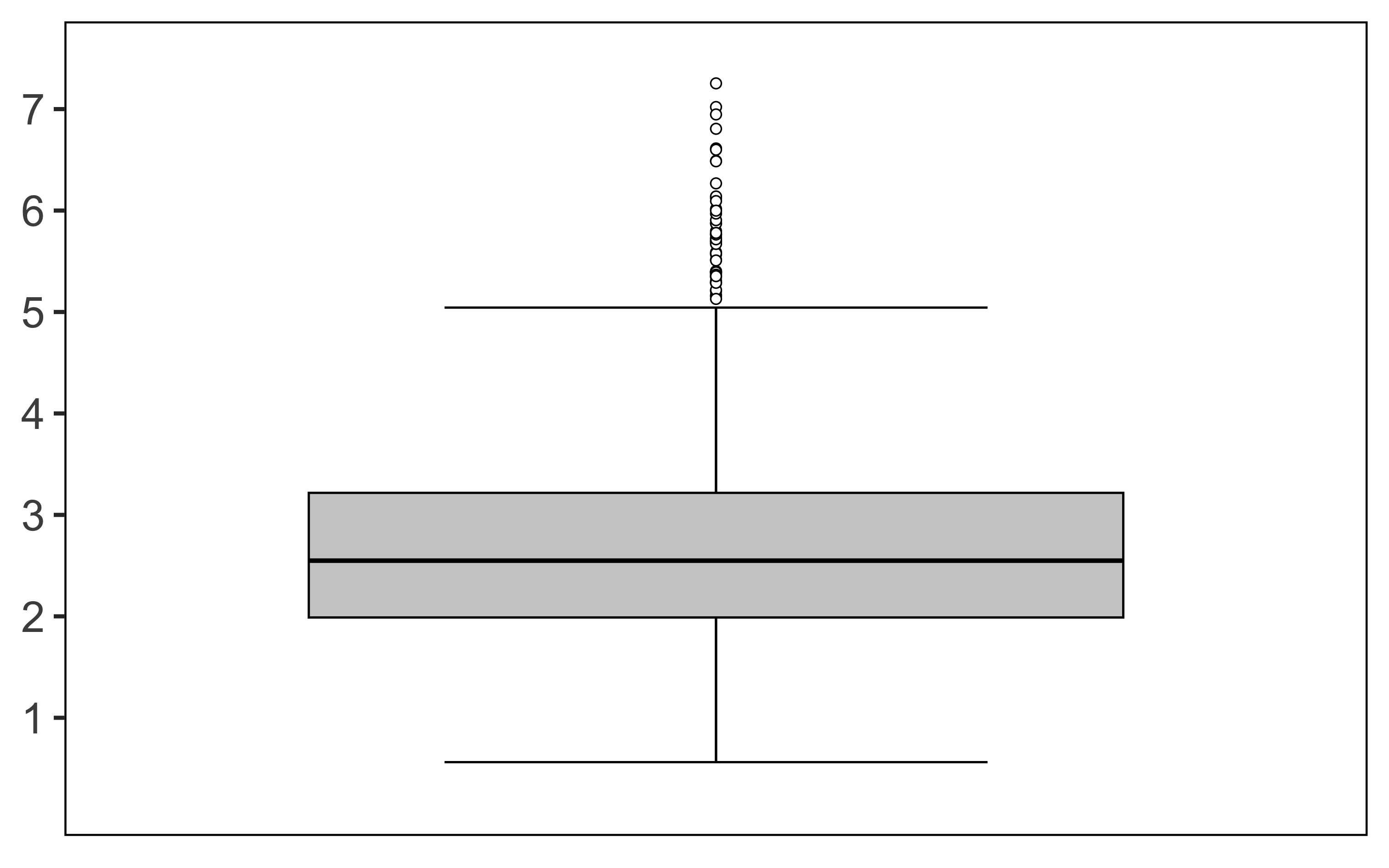}
    \caption{Boxplot of estimated GMDs per county.}
    \label{fig:boxplot_GMD}
\end{figure}

We begin by evaluating whether the Gini mean difference of the income-to-hours ratio (\texttt{INCTOT}/\texttt{UHRSWORK}) is homogeneous across counties. This measure serves as a proxy for effective hourly pay and may reveal whether income disparities vary significantly across regions. 
The estimated GMDs range from 0.56 to 7.25. Figure~\ref{fig:boxplot_GMD} displays their distribution across counties.

\noindent At this point, we apply the resampling methodology described in Subsection \ref{sect:6.2}.   In Table \ref{tab:Resampling_gini}, we report results for values of $L$ ranging from 5 to 40 in increments of 5 and for $\tilde{n}\in\{10,20,30,40,50\}$. 



\begin{table}[h]
    \centering
    \resizebox{\linewidth}{!}{%
    \begin{tabular}{l c c c c c c c c c}
      \toprule
      $\tilde{n}$ & Test & $L=5$ & $L=10$ & $L=15$ & $L=20$ & $L=25$ & $L=30$ & $L=35$ & $L=40$ \\
      \midrule
      \multirow{2}{*}{10}
        & N & $9.93 \cdot 10^{-3}$  & $4.24 \cdot 10^{-5}$ & $9.27 \cdot 10^{-4}$  & $1.76 \cdot 10^{-12}$ & $8.76 \cdot 10^{-20}$ & $6.13 \cdot 10^{-20}$ & $1.01 \cdot 10^{-8}$ & $8.24 \cdot 10^{-9}$ \\
        & LB & 0.03  & 0  & 0.40 & 0 & 0 & 0 & 0 & 0 \\
      \midrule
      \multirow{2}{*}{20}
        & N & $3.41 \cdot 10^{-3}$ & $6.13 \cdot 10^{-21}$ & $9.35 \cdot 10^{-21}$ & $8.20 \cdot 10^{-18}$ & $1.44 \cdot 10^{-16}$ & $9.30 \cdot 10^{-28}$ & $1.32 \cdot 10^{-27}$ & $2.67 \cdot 10^{-32}$ \\
        & LB & 0.22 & 0 & 0 & 0 & 0 & 0 & 0 & 0 \\
      \midrule
      \multirow{2}{*}{30}
        & N & $2.47 \cdot 10^{-25}$ & $1.25 \cdot 10^{-14}$ & $3.29 \cdot 10^{-20}$ & $2.97 \cdot 10^{-41}$ & $3.33 \cdot 10^{-37}$ & $1.56 \cdot 10^{-38}$ & $8.07 \cdot 10^{-26}$ & $5.96 \cdot 10^{-42}$ \\
        & LB & 0 & 0 & 0 & 0 & 0 & 0 & 0 & 0 \\
      \midrule
      \multirow{2}{*}{40}
        & N & $2.96 \cdot 10^{-26}$ & $4.61 \cdot 10^{-43}$ & $3.12 \cdot 10^{-42}$ & $8.43 \cdot 10^{-41}$ & $3.03 \cdot 10^{-25}$ & $8.57 \cdot 10^{-46}$ & $5.64 \cdot 10^{-37}$ & $2.89 \cdot 10^{-49}$ \\
        & LB & 0 & 0 & 0 & 0 & 0 & 0 & 0 & 0 \\
      \midrule
      \multirow{2}{*}{50}
        & N & $4.56 \cdot 10^{-30}$  & $3.41 \cdot 10^{-28}$ & $5.67 \cdot 10^{-45}$ & $7.69 \cdot 10^{-38}$ & $7.69 \cdot 10^{-43}$ & $1.86 \cdot 10^{-40}$ & $5.30 \cdot 10^{-44}$ & $4.78 \cdot 10^{-50}$ \\
        & LB & 0 & 0 & 0 & 0 & 0 & 0 & 0 & 0 \\
      \bottomrule
    \end{tabular}%
    }
    \caption{Adjusted $p$-values for balanced group sizes $\tilde{n}$ using $L$ random samples, based on the normal approximation (N) and the linear bootstrap approximation (LB).}
    \label{tab:Resampling_gini}
\end{table}
\noindent Table~\ref{tab:Resampling_gini} shows that for \( L \geq 20 \) and \( \tilde{n} \geq 10 \), the test consistently rejects the null hypothesis, with associated \( p \)-values being extremely small. Moreover, for \( \tilde{n} \geq 30 \), the decision is stable across all considered values of \( L \). This observation is consistent with the results in Section~\ref{section:simulation}: under alternatives, the power increases with \( n_0 \), while under the null, the level is only slightly affected by the choice of \( n_0 \). Altogether, the results provide strong evidence against the null hypothesis and support the practical adequacy of using \( L = 30 \) as an upper bound, as originally proposed in \cite{Cuesta2010}.

\subsection*{Spearman's correlation between working hours and earned income}\label{subsec:spearman_correlation}

\begin{figure}[ht]
    \centering
    \includegraphics[width=0.7\linewidth]{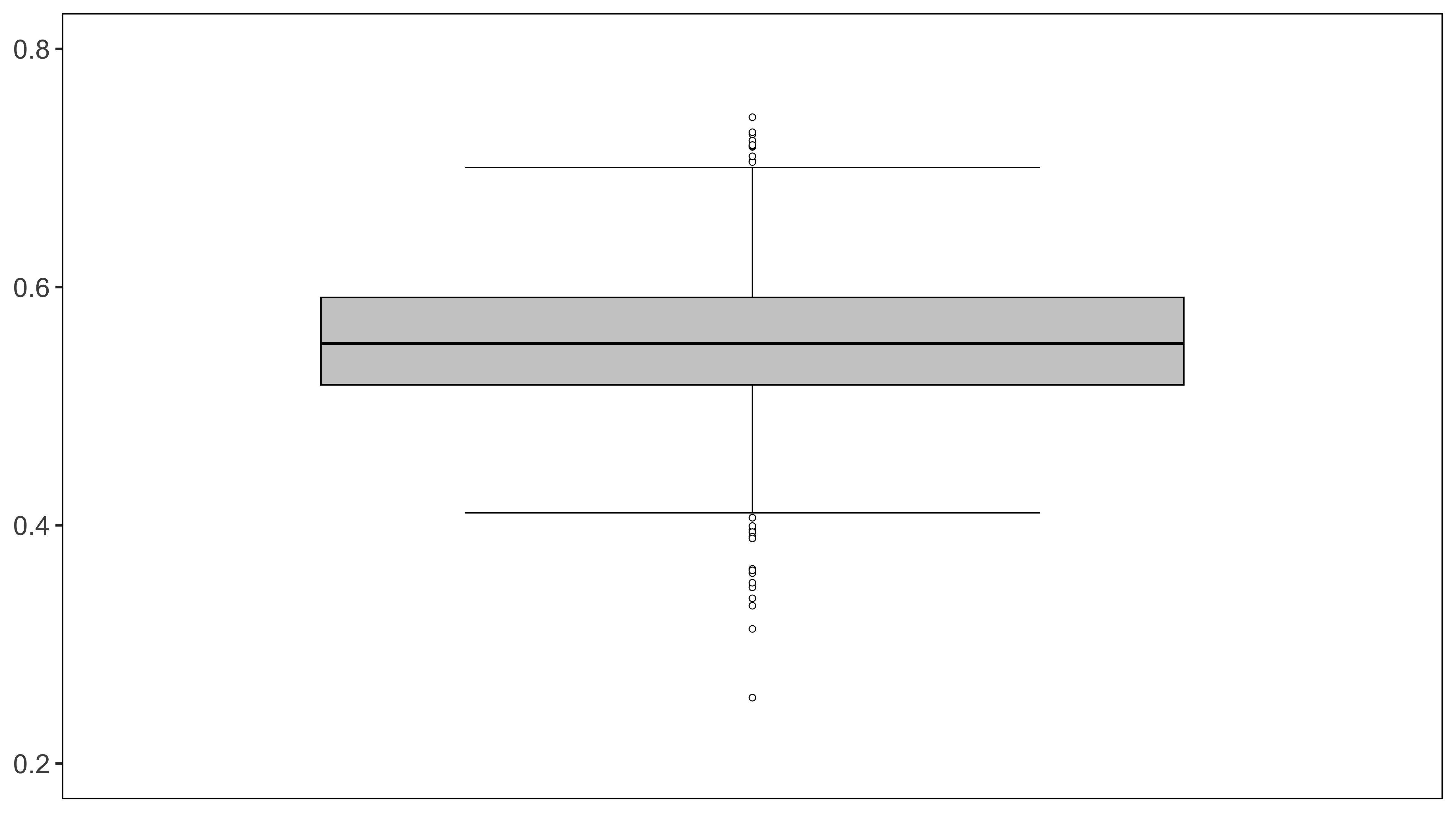}
    \caption{Boxplot of approximated Speraman's rhos per county.}
    \label{fig:boxplot_sprh}
\end{figure}

We analyze $\rho_S$ between the number of usual hours worked per week and earned income to assess the homogeneity of the relationship across all counties. The estimated Spearman’s rho are all positive, ranging from 0.26 to 0.74, indicating that in every county the two variables exhibit a positive association. The distribution of estimated $\rho_S$ is displayed in Figure \ref{fig:boxplot_sprh}.

As done in the previous case, we now perform the same resampling study. The resulting $p$-values are presented in Table~\ref{tab:Resampling_sprho}.

\begin{table}[ht]
    \centering
    \resizebox{\linewidth}{!}{%
    \begin{tabular}{l c c c c c c c c c}
      \toprule
      $\tilde{n}$ & Test & $L=5$ & $L=10$ & $L=15$ & $L=20$ & $L=25$ & $L=30$ & $L=35$ & $L=40$ \\
      \midrule
      \multirow{2}{*}{10}
        & N & 0.60  & 0.75 & 0.74  & 0.36 & 0.51 & 1 & 0.46 & 0.37 \\
        & LB & 0.90  & 1  & 1 & 1 & 1 & 1 & 1 & 1 \\
      \midrule
      \multirow{2}{*}{20}
        & N & 0.02 & $5.20 \cdot 10^{-3}$ & $2.17 \cdot 10^{-3}$ & $3.33 \cdot 10^{-4}$ & $4.9 \cdot 10^{-3}$ & 0.02 & $1.41 \cdot 10^{-3}$ & $1.74 \cdot 10^{-3}$ \\
        & LB & 0.40 & 0.19 & 0.12 & 0 & 0 & 0.78 & 0.25 & 0 \\
      \midrule
      \multirow{2}{*}{30}
        & N & $8.11 \cdot 10^{-3}$ & $1.85 \cdot 10^{-4}$ & $1.68 \cdot 10^{-4}$ & $1.83 \cdot 10^{-7}$ & $1.96 \cdot 10^{-4}$ & $1.06 \cdot 10^{-4}$ & $7.88 \cdot 10^{-5}$ & $8.41 \cdot 10^{-5}$ \\
        & LB & 0.01 & 0.13 & 0 & 0 & 0 & 0 & 0 & 0 \\
      \midrule
      \multirow{2}{*}{40}
        & N & $1.02 \cdot 10^{-6}$ & $4.76 \cdot 10^{-7}$ & $8.30 \cdot 10^{-6}$ & $9.57 \cdot 10^{-7}$ & $3.26 \cdot 10^{-6}$ & $5.87 \cdot 10^{-8}$ & $1.48 \cdot 10^{-7}$ & $5.640 \cdot 10^{-37}$ \\
        & LB & 0.18 & 0 & 0 & 0 & 0 & 0 & 0 & 0 \\
      \midrule
      \multirow{2}{*}{50}
        & N & $3.44 \cdot 10^{-8}$  & $4.08 \cdot 10^{-8}$ & $6.38 \cdot 10^{-8}$ & $1.72 \cdot 10^{-9}$ & $7.51 \cdot 10^{-10}$ & $3.65 \cdot 10^{-9}$ & $7.58 \cdot 10^{-9}$ & $1.09 \cdot 10^{-8}$ \\
        & LB & 0 & 0 & 0 & 0 & 0 & 0 & 0 & 0 \\
      \bottomrule
    \end{tabular}%
    }
    \caption{Adjusted $p$-values for balanced group sizes $\tilde{n}$ using $L$ random samples, based on the normal approximation (N) and the linear bootstrap approximation (LB).}
    \label{tab:Resampling_sprho}
\end{table}
Looking at this table, we see that for \( \tilde{n} \geq 30 \) and \( L \geq 15 \), the test consistently rejects the null, with low \( p \)-values. Moreover, for \( \tilde{n} = 50 \), $H_0$ is rejected for every considered value of \( L \). In contrast, for \( \tilde{n} = 10 \), the test does not reject the null, even for large values of \( L \). These findings are in line with those discussed for the Gini case and more clearly highlight the influence of \( \tilde{n} \) on the test decision, as previously noted. Taken together, the results provide strong evidence against the null in this setting.

\section{Conclusions and future work}
\label{section:conclusions}

In this work, a general, nonparametric test based on $U$–statistics was developed for assessing the equality of parameters across a large number of populations. Under mild moment and sample‐size comparability conditions, the proposed test statistic converges to a normal distribution under the null hypothesis. In addition, a ratio‐consistent estimator of its variance was provided, which, together with the statistic’s convergence, enables asymptotically exact inference without parametric assumptions. A linear bootstrap approximation to improve the performance for small to moderate $k$, and a random sampling strategy to reduce computational burden, were also implemented. Simulation studies for the Gini mean difference and Spearman’s rho were conducted to examine the behavior of both the asymptotic and bootstrap versions in terms of Type I error control and power, across different group sizes and alternatives. A real-data application was also presented to illustrate the procedure.   

A number of extensions are possible. In many real‐world settings the observations within each group may exhibit temporal or spatial dependence (for example, in time‐series or panel data), so extending the theory to accommodate dependent data is an important avenue for future work. Likewise, it is often of interest to compare vector‐valued or multivariate functionals, such as covariance matrices or other high‐dimensional summaries, across many populations. In this regard, generalizing the current framework to parameters of arbitrary dimension is therefore another line for future research.

\section{Proofs}
\label{section: proofs}
Throughout this section, \(M, M_1, M_2, \dots\) denote generic positive constants that are not important and may vary from one step to another, or even between different occurrences within the same line of an inequality. A sequence of random variables is $o_{\Prob}(1)$ if it converges to 0, and is $O_{\Prob}(1)$ if it is bounded in probability.

\begin{proof}[Proof of Lemma \ref{lemma:consistency}]
    Note that
    \[T_k=\left(\frac{1}{k}-\frac{1}{k^2}\right) \sum_{i=1}^k \widehat{\theta^2_i}-\Bigg(\frac{1}{k} \sum_{i=1}^k \widehat{\theta_i} \Bigg)^2+\frac{1}{k^2} \sum_{i=1}^k \widehat{\theta_i}^2.\]
We know that $\E(\widehat{\theta}_i)=\theta_i$, $\E(\widehat{\theta^2_i})=\theta^2_i$, $\E\Big(\widehat{\theta_i}^2\Big)=\theta_i^2+\V(\widehat{\theta_i})$  for $\forall i$.
To prove the first part of the result it suffices to show 
\begin{equation}
    \frac{1}{k} \sum_{i=1}^k \widehat{\theta_i}-\frac{1}{k} \sum_{i=1}^k \theta_i \tas 0,
    \label{eq:consistency1}
\end{equation}
\begin{equation}
    \frac{1}{k} \sum_{i=1}^k \widehat{\theta^2_i}-\frac{1}{k} \sum_{i=1}^k \theta^2_i \tas 0,
    \label{eq:consistency2}
\end{equation}
and
\begin{equation}
    \frac{1}{k^2} \sum_{i=1}^k \widehat{\theta_i}^2 \tp 0.
    \label{eq:consistency3}
\end{equation}
Firstly we prove \eqref{eq:consistency1} and \eqref{eq:consistency2}. As the random variables $\widehat{\theta_1},\ldots, \widehat{\theta_k}$ (also $\widehat{\theta^2_1},\ldots, \widehat{\theta^2_k}$)  are independent but not identically distributed we apply Kolmogorov's Strong Law of Large Numbers (see Theorem D in Chapter 1.8 of \cite{serfling1980approximation}). To do so, it suffices to show that 
$\V(\widehat{\theta_i}), \V(\widehat{\theta^2_i}) \leq M, \, \forall i.$
 From \eqref{eq:sigmaorder} and \eqref{eq:varsigma}, we have
\begin{equation}
    \V(\widehat{\theta_i})\leq M_1 \sigma_m^2= M_1 \V\{h(X_{i1},\ldots, X_{im})\} \leq M_1 \E\{h^2(X_{i1},\ldots, X_{im})\}.
    \label{eq:variancetheta1bounded}
\end{equation}
An analogous equation holds for 
\begin{align}
\V(\widehat{\theta^2_i}) &\leq M_1 \sigma_{2m}^2 = M_1 \V\{H(X_{i1},\ldots,X_{i2m})\} \notag \\
&\leq M_1 \E\{H^2(X_{i1},\ldots,X_{i2m})\} \notag \\
&\leq M_1 \E\{ h^2(X_{i1},\ldots,X_{im}) h^2(X_{i(m+1)},\ldots,X_{i2m})\} \notag \\
&= M_1\E^2\{ h^2(X_{i1},\ldots, X_{im}) \},
\label{eq:variancetheta2bound}
\end{align}
where the last inequality follows from the independence of $\{X_{im}\ldots,X_{im}\}$ and $\{X_{i(m+1)},\ldots,X_{i2m}\}$.

Therefore, since we have \eqref{ass:boundenessh2}, we obtain  \eqref{eq:consistency1} and \eqref{eq:consistency2}, which in turn implies convergence in probability. Finally, we proceed to prove \eqref{eq:consistency3}. From \eqref{ass:boundenessh2} and \eqref{eq:variancetheta1bounded}, and observing that
\begin{equation}
    \theta_i=\E\{h(X_{i1},\ldots,X_{im})\}\leq \E^{1/2}\{h^2(X_{i1},\ldots,X_{im})\},
    \label{eq:thetaacotado}
\end{equation}
it follows that $\E\left(\widehat{\theta_i}^2\right)=\theta_i^2+\V(\widehat{\theta_i})<M, \, \forall i$, we get  $(1/k)\sum_{i=1}^k\widehat{\theta_i}^2=O_{\Prob}(1)$ and then \eqref{eq:consistency3} is fulfilled.

For the second part, we only need to show that the convergence in \eqref{eq:consistency3} is almost sure.
Note that
\[\frac{1}{k^2} \sum_{i=1}^k \widehat{\theta_i}^2=\frac{1}{k}\left\{\frac{1}{k}\sum_{i=1}^k \widehat{\theta_i}^2- \frac{1}{k}\sum_{i=1}^k \E\left(\widehat{\theta_i}^2\right)\right\}+ \frac{1}{k^2} \sum_{i=1}^k\E\left(\widehat{\theta_i}^2\right).\]
From \eqref{eq:variancetheta1bounded}, we know that $E\left(\widehat{\theta_i}^2\right)\leq M,\, \forall i$, so to prove the result, it suffices to show that
\[\frac{1}{k}\sum_{i=1}^k \widehat{\theta_i}^2- \frac{1}{k}\sum_{i=1}^k \E\left(\widehat{\theta_i}^2\right) \tas 0.\]
To do so, we apply Corollary 3.9 of \cite{white2014asymptotic}, where we only need to check that $\E(|\widehat{\theta_i}|^{2+\delta}) \leq M, \, \forall i$, for some $\delta>0$ and all $\forall i$. Using the $c_r$ inequality, we get

\begin{align}
\E(|\widehat{\theta_i}|^{2+\delta})&= \frac{1}{\{n_i(m)\}^{2+\delta}}\E\left\{\Big| \sum_{j_1 \neq \cdots \neq j_m} h(X_{ij_1}, \ldots, X_{ij_m})\Big|^{2+\delta}\right\} \nonumber \\ 
     &\leq \frac{1}{n(m)} \sum_{j_1 \neq \cdots \neq j_m} \E\left\{|h(X_{ij_1}, \ldots, X_{ij_m})|^{2+\delta}\right\} \nonumber\\&= \E\left\{|h(X_{ij_1}, \ldots, X_{ij_m})|^{2+\delta}\right\}.
\label{eq:expwidthetadelta}
\end{align}
Then, as $\E\{ |h(X_{ij_1}, \ldots, X_{ij_m})|^{2+\delta}\} < M$ for $\forall i$, the result is proven.
\end{proof}

\begin{proof}[Proof of Lemma \ref{lemma:Tkdecomp}]
To obtain the result, we follow a similar approach to that described in Section 3 of \cite{AkritasPapadatos2004}.
In this case, the projection is applied by conditioning to the samples \( \bm{X}_1, \ldots, \bm{X}_k \), rather than to the individual random variables. By applying the Decomposition Lemma of \cite{efron1981jackknife} we get
\begin{equation}
    T_k = D_k + \sum_{r=1}^k \left\{ \E(T_k | \bm{X}_r) - D_k \right\} + \sum_{r < s}^k \left\{ \E(T_k | \bm{X}_r, \bm{X}_s) - \E(T_k | \bm{X}_r) - \E(T_k | \bm{X}_s) + D_k \right\}.
    \label{eq:Hoeffdingapplied}
\end{equation}

Computing the conditional expectations, we find
\begin{align*}
    \E(T_k | \bm{X}_r)&= \left(\frac{1}{k} - \frac{1}{k^2}\right) \Big(\sum_{i\neq r} \theta_i^2 + \widehat{\theta_r^2}\Big)-\frac{1}{k^2}\Big( \sum_{i\neq j \neq r} \theta_i \theta_j +2 \widehat{\theta_r}\sum_{i\neq r} \theta_i \Big) \\ 
    &= \underbrace{\left(\frac{1}{k} - \frac{1}{k^2}\right) \sum_{i=1}^k \theta_i^2 -\frac{1}{k^2} \sum_{i\neq j } \theta_i \theta_j}_{=D_k} +\left(\frac{1}{k} - \frac{1}{k^2}\right) (\widehat{\theta_r^2} - \theta_r^2) - \frac{2}{k^2}(\widehat{\theta_r}-\theta_r)\Big(\sum_{i=1}^k \theta_i -\theta_r\Big) \\
    &=D_k + \left(\frac{1}{k} - \frac{1}{k^2}\right)(\widehat{\theta^2_r} - \theta_r^2) - \frac{2}{k}\overline{\theta}(\widehat{\theta}_r - \theta_r) + \frac{2}{k^2}\theta_r(\widehat{\theta}_r - \theta_r) \\
    & =D_k+\frac{1}{k}L_{rk}-\frac{1}{k^2}(\widehat{\theta^2_r} - \theta_r^2)+\frac{2}{k^2}\theta_r(\widehat{\theta}_r - \theta_r),
\end{align*}
and  for $r\neq s$
\begin{align*}
    \E(T_k | \bm{X}_r, \bm{X}_s) &= \left(\frac{1}{k} - \frac{1}{k^2}\right) \Big(\sum_{i\neq r,s} \theta_i^2 + \widehat{\theta_r^2}+ \widehat{\theta_s^2}\Big)-\frac{1}{k^2}\Big( \sum_{i\neq j \neq r,s} \theta_i \theta_j +2 \widehat{\theta_r}\sum_{i\neq r} +2 \widehat{\theta_s}\sum_{i\neq s} \theta_i + 2 \widehat{\theta_r}  \widehat{\theta_s} \Big) \\
    &= \underbrace{\left(\frac{1}{k} - \frac{1}{k^2}\right) \sum_{i=1}^k \theta_i^2 -\frac{1}{k^2} \sum_{i\neq j } \theta_i \theta_j}_{=D_k}  +\left(\frac{1}{k} - \frac{1}{k^2}\right) (\widehat{\theta_r^2} - \theta_r^2) + \left(\frac{1}{k} - \frac{1}{k^2}\right) (\widehat{\theta_s^2} - \theta_s^2) \\
    &\quad - \frac{2}{k^2}(\widehat{\theta_r}-\theta_r)\sum_{i\neq r,s} \theta_i - \frac{2}{k^2}(\widehat{\theta_s}-\theta_s)\sum_{i\neq r,s} \theta_i - \frac{2}{k^2}(\widehat{\theta_r} \widehat{\theta_s}-\theta_r\theta_s) \\
    &= D_k + \underbrace{\left(\frac{1}{k} - \frac{1}{k^2}\right) (\widehat{\theta_r^2} - \theta_r^2) - \frac{2}{k}\overline{\theta}(\widehat{\theta}_r - \theta_r)+ \frac{2}{k^2}\theta_r(\widehat{\theta}_r - \theta_r)}_{ E(T_k | \bm{X}_r)-D_k} \\
    &\quad + \underbrace{\left(\frac{1}{k} - \frac{1}{k^2}\right)(\widehat{\theta_s^2} - \theta_s^2) - \frac{2}{k} \overline{\theta} (\widehat{\theta}_s - \theta_s) + \frac{2}{k^2} \theta_s (\widehat{\theta}_s - \theta_s)}_{ E(T_k | \bm{X}_s)-D_k} \\
    &\quad + \frac{2}{k^2} \theta_r (\widehat{\theta}_s - \theta_s) + \frac{2}{k^2} \theta_s (\widehat{\theta}_r - \theta_r) - \frac{2}{k^2}(\widehat{\theta_r} \widehat{\theta_s}-\theta_r\theta_s) \\
    &= \E(T_k | \bm{X}_r) + \E(T_k | \bm{X}_s) -D_k - \frac{2}{k^2}(\widehat{\theta}_r - \theta_r)(\widehat{\theta}_s - \theta_s).
\end{align*}
Substituting these expressions into \eqref{eq:Hoeffdingapplied}, the result follows.
\end{proof}

\begin{remark}
Note that \( L_{ik} \) can be expressed as the sum of \( \mathbb{L}_{i1} \) and \( \mathbb{L}_{ik2} \), where
\begin{equation}
    \begin{aligned}
        \mathbb{L}_{i1} &= \widehat{\theta_i^2} - \theta_i^2 - 2\theta_i(\widehat{\theta_i} - \theta_i), \\
        \mathbb{L}_{ik2} &= 2(\theta_i - \overline{\theta})(\widehat{\theta_i} - \theta_i).
    \end{aligned}
    \label{eq:L_terms}
\end{equation}
The term \( \mathbb{L}_{i1} \) corresponds to a degree-$2m$ \( U \)-statistic with a symmetric kernel \( F \),
\begin{align*}
F(x_1,\ldots,x_{2m})=\theta_i^2 + H(x_1,\ldots,x_{2m})-2\theta_i\binom{2m}{m}^{-1}\sum_{1\leq j_1<\cdots<j_{m}\leq 2m}h(x_{j_1},\ldots,x_{j_m}),
\end{align*}
where \( H \) is defined in \eqref{eq:H}. We have that 
\begin{align*}
    \E\{F(X_{i1},\ldots,X_{i2m})|X_{i1}\} &= \theta_i^2 + \theta_i h_1(X_{i1})  \\&\quad - 2\theta_i\binom{2m}{m}^{-1} \left[ \binom{2m-1}{m-1} h_1(X_{i1}) + \left\{ \binom{2m}{m} - \binom{2m-1}{m-1} \right\} \theta_i \right] \\ 
    &= 0,
\end{align*}
since $\binom{2m}{m}=2\binom{2m-1}{m-1}$. Therefore \( \mathbb{L}_{i1} \) is degenerated. Now, we compute the second order variance, i.e. $\V[\E\{F(X_{i1},\ldots,X_{i2m})|X_{i1},X_{i2}\}]$. First, we calculate the conditional expectation
\begin{equation}
\begin{aligned}
    F_2(x,y)&=\E\{F(X_{i1},\ldots,X_{i2m})|X_{i1}=x,X_{i2}=y\}\\&=\theta_i^2+H_2(x,y) -2\theta_i\binom{2m}{m}^{-1} \sum_{1\leq j_1<\cdots<j_{m}\leq 2m} E\{h(X_{ij_1},\ldots,X_{ij_m})|X_{i1}=x,X_{i2}=y\},
\end{aligned}
\label{eq:F2general}
\end{equation}
 where $H_2(x,y)=\E\{H(X_{i1},\ldots,X_{i2m})|X_{i1}=x,X_{i2}=y\}$. Routine calculations show that
\begin{equation}
    H_2(x,y)=2\binom{2m}{m}^{-1}\left\{\binom{2m-2}{m-2}\theta_i h_2(x,y)+\binom{2m-2}{m-1}h_1(x)h_1(y)\right\}.
    \label{eq:Hc2}
\end{equation}
On the other hand,
\begin{equation}
\begin{array}{c}
\displaystyle
     \sum_{1 \leq j_1 < \cdots < j_m \leq 2m} \E\{h(X_{ij_1}, \ldots, X_{ij_m}) | X_{i1}=x, X_{i2}=y\}=   \\
     \displaystyle
     \binom{2m-2}{m-1} \left\{ h_1(x) + h_1(y) \right\} 
 +\binom{2m-2}{m-2} h_2(x, y) 
+ \left\{\binom{2m}{m}-2\binom{2m-2}{m-1}-\binom{2m-2}{m-2}\right\}\theta_i. 
\end{array}
\label{eq:sumExpect}
\end{equation}
Substituting \eqref{eq:Hc2} and \eqref{eq:sumExpect} into \eqref{eq:F2general}, it follows that
\[F_2(x,y) =\frac{m}{2m-1}\left[h_1(x)h_1(y)-\theta_i\left\{h_1(x)+h_1(y)\right\}+\theta_i^2\right],\]
and noting that,
\begin{equation*}
   \begin{array}{c}
   \displaystyle
        \C\Big\{h_1(X_{i1})h_1(X_{i2}),h_1(X_{i1})+h_1(X_{i2})\Big\}= \\
        \displaystyle
         2 \C\Big\{h_1(X_{i1})h_1(X_{i2}),h_1(X_{i1})\Big\}
    = 2 \E\{h_1(X_{i2})\}\V\{h_1(X_{i1})\} = 2 \theta_i \sigma_{i1}^2
   \end{array} 
\end{equation*}
we get

\[
\V\left\{F_2(X_{i1},X_{i2})\right\} = \frac{m^2}{(2m-1)^2} \sigma_{i1}^4.
\]
Since $\mathbb{L}_{i1}$ is degenerated, from \eqref{eq:varsigma}, \eqref{eq:compsamp} and the previous equation, it follows that
\begin{equation}
n_0^2 \V(\mathbb{L}_{i1})\geq M\sigma_{i1}^4.
    \label{eq:lowerLi1}
\end{equation}
On the other hand, using \eqref{eq:sigmaorder},
\eqref{eq:varsigma}  and \eqref{eq:compsamp}, we obtain  
\begin{equation}
   n_0^2 \V(\mathbb{L}_{i1}) \leq M \V\{F(X_{i1},\ldots,X_{i2m})\}.
   \label{eq:boundLi1prev}
\end{equation}
By applying $c_r$ inequality we get
\begin{align*}
    \V\{F(X_{i1},\ldots,X_{i2m})\} &\leq  \E\{F^2(X_{i1},\ldots,X_{i2m})\} \\
    &\leq M\Big[\E\{H^2(X_{i1},\ldots,X_{i2m})\}+\theta_i^2 \E\{h^2(X_{i1},\ldots,X_{im})\}\Big],
\end{align*} Furthermore, for the monotonicity of moments, we note that
\begin{equation}
\theta_i\leq |\theta_i|=|\E\{h(X_{i1},\ldots,X_{im})\}| \leq \E\{|h(X_{i1},\ldots,X_{im})|\} \leq \E^{1/2}\{h^2(X_{i1},\ldots,X_{im})\}.
\label{eq:boundthetai}
\end{equation}
In addition, from \eqref{eq:variancetheta2bound}, we have  
\[
\E\{H^2(X_{i1},\ldots,X_{i2m})\} \leq M \E^2\{h^2(X_{i1},\ldots,X_{im})\}.
\]  
Substituting all the previous bounds into \eqref{eq:boundLi1prev}, we obtain  
\begin{equation}
n_0^2 \V(\mathbb{L}_{i1}) \leq M \E^2\{h^2(X_{i1},\ldots,X_{im})\}.
    \label{eq:upperLi1}
\end{equation}

The term $\mathbb{L}_{ik2}$ is a degree-2 $U$-statistic whose kernel is given by
\[f(x_1,\ldots,x_m)=2(\theta_i-\overline{\theta})\{h(x_1,\ldots,x_m)-\theta_i\},\]
and then 
\begin{equation}
    \V\{f_1(X_{i1})\}=\V[\E\{f(X_1,\ldots,X_m)|X_1\}]=4(\theta_i-\overline{\theta})^2\sigma_{i1}^2.
    \label{eq:varf}
\end{equation}
This implies that, in general, \( \mathbb{L}_{ik2} \) is a non-degenerated \( U \)-statistic, and then,  for analogous reasons to \( \mathbb{L}_{i1} \), we have that
\begin{equation}
n_0 \V(\mathbb{L}_{ik2})\geq M(\theta_i-\overline{\theta})^2\sigma_{i1}^2.
    \label{eq:lowerL2ik}
\end{equation}
We also have an upper bound on $\V(\mathbb{L}_{ik2})$. From \eqref{eq:sigmaorder}, \eqref{eq:varsigma} and \eqref{eq:compsamp}, we get
\begin{equation}
    n_0 \V(\mathbb{L}_{ik2}) \leq M \V\{f(X_{i1},\ldots,X_{im})\}\leq  M(\theta_i-\overline{\theta})^2 \E\{h^2(X_{i1},\ldots,X_{im})\},
    \label{eq:upperLik2theta}
\end{equation} 
then, by $c_r$ inequality  and \eqref{eq:boundthetai}, it follows that 
\begin{align*}
    \V\{f(X_{i1},\ldots,X_{im})\} &= 4(\theta_i-\overline{\theta})^2 \V\{h(X_{i1},\ldots,X_{im})\} \\
    &\leq 8 (\theta_i^2 + \overline{\theta}^2) \E\{h^2(X_{i1},\ldots,X_{im})\}\\ 
    &\leq M \E^2\{h^2(X_{i1},\ldots,X_{im})\},
\end{align*}
and hence
\begin{equation}
    n_0 \V(\mathbb{L}_{ik2}) \leq M \E^2\{h^2(X_{i1},\ldots,X_{im})\}.
    \label{eq:upperLik2}
\end{equation}
\label{remark}
\end{remark}

Before proving Lemma~\ref{lemma:variancegeneral}, we first present some preliminary results. Throughout this work, by convention, we set the binomial coefficient as
\begin{equation}
    \binom{a}{b} = 0 \quad \text{if} \quad b > a.
    \label{convention}
\end{equation}
In the following proof, we will make use of the identities provided in Section~\ref{section:binomials} of the appendix.
\begin{lemma}
    Assume that \eqref{eq:independent_sample} holds and that $\E\{h^2(X_{1},\ldots,X_{m})\}< \infty, \; \forall i$. Let  \
    \[
    S_m=\binom{n}{2m}^{-1}\binom{2m}{m}^{-1} \sum_{c=1}^{m} \binom{n-m}{2m-c}\binom{m}{c} S_{mc}, \qquad S_{mc}=\sum_{j=0}^{c} \binom{c}{j}  \binom{2m-c}{m-j} \left(\theta\sigma_j^2+\theta\sigma_{c-j}^2\right).
    \]
    Then, \( S_m = \V(2\widehat{\theta})/2 \).
    \label{lemma:Vandermondeapp1}
\end{lemma}
\begin{proof}
First, from \eqref{eq:nkkj} we note that
\begin{equation}
    \binom{n}{2m}^{-1}\binom{2m}{m}^{-1}=\binom{n}{m}^{-1}\binom{n-m}{m}^{-1}.
    \label{eq:2mm-1}
\end{equation}
In addition, we have that
\begin{align}
S_{mc}=
\sum_{j=0}^{c} \binom{c}{j}  \binom{2m-c}{m-j} \theta\sigma_j^2+ \sum_{s=0}^{c} \binom{c}{c-s}  \binom{2m-c}{m-(c-s)} \theta\sigma_{s}^2.
\label{eq:combinesums}
\end{align}
Then, from the symmetry of the binomial coefficients and noting that the term for $j=0$ equals 0 since $\sigma_0=0$, we get that
\begin{align*}
S_{mc}=\sum_{j=1}^{c} \binom{c}{j}  \binom{2m-c}{m-j} \theta\sigma_j^2+ \sum_{s=1}^{c} \binom{c}{s}  \binom{2m-c}{m-s} \theta\sigma_{s}^2=2 \sum_{j=0}^{c} \binom{c}{j}  \binom{2m-c}{m-j} \theta\sigma_j^2.
\end{align*}
Hence, using the last equation we can write
\begin{align*}
    S_m&=2\theta\binom{n}{m}^{-1}\binom{n-m}{m}^{-1} \sum_{c=1}^{m} \binom{n-m}{2m-c}\binom{m}{c} \sum_{j=1}^{c} \binom{c}{j}  \binom{2m-c}{m-j} \sigma_j^2 \\
    &= 2\theta\binom{n}{m}^{-1}\binom{n-m}{m}^{-1} \sum_{j=1}^{m} \sigma_j^2\sum_{c=j}^{m} \binom{m}{c}  \binom{c}{j}  \binom{n-m}{2m-c}\binom{2m-c}{m-j},
\end{align*}
where we have changed the order of summations. Now applying property \eqref{eq:nkkj} twice, we get
\[S_m = 2\theta\binom{n}{m}^{-1}\binom{n-m}{m}^{-1} \sum_{j=1}^{m} \binom{m}{j} \binom{n-m}{m-j} \sigma_j^2 \sum_{c=j}^{m} \binom{m-j}{c-j}  \binom{n-2m+j}{m-(c-j)} \]
Now we apply Vandermonde's identity \eqref{eq:vandermonde} yielding
\begin{align}
     \sum_{c=j}^{m} \binom{m-j}{c-j}  \binom{n-2m+j}{m-(c-j)}   =\sum_{s=0}^{m-j} \binom{m-j}{s}  \binom{n-2m+j}{m-s}=\sum_{s=0}^{m} \binom{m-j}{s}  \binom{n-2m+j}{m-s}=\binom{n-m}{m},
\end{align}
where we have extended the sum to $m$ because $\binom{m-j}{s}=0$ if $s>m-j$. Then, recalling the definition of $\V(\widehat{\theta})$ in \eqref{eq:varsigma}, we obtain the desired result.

\end{proof}

\begin{lemma}
Assume that \eqref{eq:independent_sample} holds and that $\E\{h^2(X_{1},\ldots,X_{m})\}< \infty$. $\C\left(\widehat{\theta^2}, \widehat{\theta}\right)$ can be expressed as  
\[
\C\left(\widehat{\theta^2}, \widehat{\theta}\right) =  \frac{\theta}{2}\V\left(2\widehat{\theta}\right) + \frac{1}{4} \Xi^3,
\]  
where  
\begin{equation}
    \Xi^3 = 4\binom{n}{2m}^{-1} \binom{2m}{m}^{-1} \sum_{c=1}^{m} \binom{n-m}{2m-c} \binom{m}{c} \sum_{j=1}^{c} \binom{c}{j} \binom{2m-c}{m-j} \widetilde{\zeta}^3_{c,j},
    \label{eq:Xi}
\end{equation}
and the term $\widetilde{\zeta}^3_{c,j}$ is defined as  
\[
\widetilde{\zeta}^3_{c,j} = \C\Big\{\widetilde{h}_c(X_{1},\ldots,X_{c}),\widetilde{h}_j(X_{1},\ldots,X_{j})\widetilde{h}_{c-j}(X_{j+1},\ldots,X_{c})\Big\},
\]  
with $\widetilde{h}$ defined in \eqref{eq:widetildeh}.
\label{lemma:Cov}
\end{lemma}

\begin{proof}
We have
\begin{align*}
    \C\left(\widehat{\theta}, \widehat{\theta^2}\right) = \frac{1}{n(m) n(2m)} \sum_{\substack{i_1 \neq \cdots \neq i_m \\ j_1 \neq \cdots \neq j_{2m}}} \C\Big\{ h(X_{i_1}, \ldots, X_{i_m}), h(X_{j_1}, \ldots, X_{j_m}) h(X_{j_{m+1}}, \ldots, X_{j_{2m}}) \Big\}.
\end{align*}

Since $h$ is symmetrically defined in its arguments we only need to compute
\[\C\Big\{ h(X_{i_1}, \ldots, X_{i_m}), h(X_{j_1}, \ldots, X_{j_m}) h(X_{j_{m+1}}, \ldots, X_{j_{2m}}) \Big\}\]
for fixed indices and multiply it by all the possible permutation of those. In the above expression, we have $3m$ indices for the random variables, of which at least $m$ must be distinct. It is evident that any covariance involving $3m$ distinct indices will be zero due to the independence of the random variables. Therefore, the non-zero terms will correspond to those in which $\#(\{i_1,\ldots,i_m\}\cap \{j_1,\ldots,j_{2m}\})=c$, where $c$ ranges from $1$ to $m$. Additionally, the value of these terms depends on the specific place where those indices are repeated in $h(X_{j_1},\ldots,X_{j_m})h(X_{j_{m+1}}, \ldots, X_{j_{2m}})$. Hence, if $c$ indices are repeated the non-zero terms, considering that all random variables are identically distributed, will be
\begin{equation}
\zeta^3_{c,j}=\C\Big\{h_c(X_1,\ldots,X_c),h_i(X_1,\ldots,X_j)h_{c-j}(X_{j+1},\ldots,X_c)\Big\}, \quad 0\leq j \leq c.
\label{eq:zeta3}
\end{equation}
Let $R_{c,j}$ be the number of ways in which $c$ indices can repeat, with $j$ of these $c$ indices appearing in the same configuration as in $\zeta^3_{c,j}$. With this notation we can write
\begin{equation}
    \C(\widehat{\theta}, \widehat{\theta^2})=\frac{1}{n(m)n(2m)}\sum_{c=1}^m\sum_{j=0}^c R_{c,j}\zeta^3_{c,j}.
    \label{eq:rcizeta3}
\end{equation}
In the following we calculate $R_{c,i}$. There are $n(m)$ ways of choosing $m$ indices from a set of $n$. Fixing one choice of those $m$, there are $\binom{m}{c} \binom{c}{j}$ ways to choose $c$ elements from those $m$ and partition them into two groups of $j$ and $c-j$ elements. Once this is done, we need to count the ways of distributing them into those two groups. This can be done in the same ways of those of partitioning the remaining elements, i.e. $\binom{n-m}{2m-c} \binom{2m-c}{m-j}$ ways. Finally  we have to take into accunt the $(m!)^2$ permutations of those elements. Thus, 
\begin{equation}
    R_{c,j} = n(m) \binom{m}{c} \binom{c}{j} \binom{n-m}{2m-c} \binom{2m-c}{m-j} (m!)^2,
    \label{eq:Rci}
\end{equation}
and hence, from \eqref{eq:Rci} and \eqref{eq:rcizeta3} we obtain
\begin{align}
\C(\widehat{\theta^2}, \widehat{\theta}) &= \binom{n}{2m}^{-1}\binom{2m}{m}^{-1} \sum_{c=1}^{m} \binom{n-m}{2m-c}\binom{m}{c}\sum_{j=0}^{c} \binom{c}{j}  \binom{2m-c}{m-j} \zeta^3_{c,j}.
\label{eq:covtheta2theta}
\end{align}
Now we write $\zeta^3_{c,j}$ in terms of the centered kernel $\widetilde{h}$, defined in \eqref{eq:widetildeh}. Hence, defining \begin{equation}
\widetilde{\zeta}^3_{c,j}=\C\Big\{\widetilde{h}_c(X_1,\ldots,X_c),\widetilde{h}_i(X_1,\ldots,X_j)\widetilde{h}_{c-j}(X_{j+1},\ldots,X_c)\Big\},
    \label{eq:zeta3tilde}
\end{equation} we can write
\begin{align*}
    \zeta^3_{c,j}=& \widetilde{\zeta}^3_{c,j} + \theta\C\Big\{h_c(X_{1},\ldots,X_{c}),h_j(X_{1},\ldots,X_{j})\Big\} +\theta\C\Big\{h_c(X_{1},\ldots,X_{c}),h_{c-j}(X_{j+1},\ldots,X_{c})\Big\}.
\end{align*}
Now, using \eqref{eq:sigmaelementscommon}, we can identify the covariances above, yielding
\begin{equation}
\zeta^3_{c,j}=\widetilde{\zeta}^3_{c,j}+\theta\sigma_j^2+\theta\sigma_{c-j}^2.
\label{eq:zeta3exp}
\end{equation}
Introducing \eqref{eq:zeta3exp} into \eqref{eq:covtheta2theta}, we get 
\[\C(\widehat{\theta^2}, \widehat{\theta})=\frac{1}{4}\Xi^3_{i}+S_m,\]
where
\begin{align}
    S_m=\binom{n}{2m}^{-1}\binom{2m}{m}^{-1} \sum_{c=1}^{m} \binom{n-m}{2m-c}\binom{m}{c}\sum_{j=0}^{c} \binom{c}{j}  \binom{2m-c}{m-j} (\theta\sigma_j^2+\theta\sigma_{c-j}^2).
\end{align}
From Lemma \ref{lemma:Vandermondeapp1}, $S_m= (\theta/2)\V\left(2\widehat{\theta}\right)$ and the result holds.
\end{proof}

\begin{proof}[Proof of Lemma \ref{lemma:variancegeneral}]

From Remark \ref{remark}, we have that
\begin{equation}
    \V(L_{ik})=\V(\mathbb{L}_{i1})+(\theta_i-\overline{\theta})^2\V\left(2\widehat{\theta_i}\right)+2\C(\mathbb{L}_{i1},\mathbb{L}_{ik2}),
    \label{eq:varLdecomp}
\end{equation}
and
\[\C(\mathbb{L}_{i1},\mathbb{L}_{ik2})=(\theta_i-\overline{\theta})\left\{2\C\left(\widehat{\theta_i},\widehat{\theta^2_i}\right)-\theta_i\V\left(2\widehat{\theta_i}\right)\right\}.\]
Now, recalling Lemma \ref{lemma:Cov}, the previous equation becomes
\[\C(\mathbb{L}_{i1},\mathbb{L}_{ik2})=\frac{1}{2}(\theta_i-\overline{\theta})\Xi^3_{i}.\]
Substituting this expression into \eqref{eq:varLdecomp}, the result follows.
\end{proof}

\begin{lemma} Assume that \eqref{eq:independent_sample},  \eqref{ass:boundenessh2} and \eqref{eq:compsamp} hold. Then, $n_0^2 \left|\Xi^3_{i}\right| \leq M$, $\forall i$,  and for some positive constant  $M$,
where $\Xi^3_{i}$ is defined in \eqref{eq:Xi}.
\label{lemma:Xibounded}
\end{lemma}
\begin{proof}

First, note that applying Cauchy-Schwarz inequality, we have that 
\[|\widetilde{\zeta}^3_{i(c,j)}| \leq \sigma_{ic} \sigma_{ij} \sigma_{i(c-j)}\leq \sigma_{im}^3,\]
where the last inequality follows from \eqref{eq:sigmaorder}. In addition, we know that $\sigma^2_{im}=\V\{h(X_{i1},\ldots,X_{im})\} \leq \E\{h^2(X_{i1},\ldots,X_{im})\}$, which, from \eqref{ass:boundenessh2}, is uniformly bounded, and hence
\[|\widetilde{\zeta}^3_{i(c,j)}| \leq M, \, \forall i.\]
Now, as $m$ is a fixed quantity, it follows that all the binomial coefficients involving $m$ and $c$ are bounded. Note also that $\widetilde{\zeta}^3_{i(1,j)}=0$. Then we get
    \begin{align*}
    \left|\Xi^3_{i}\right| &\leq M \binom{n_i}{2m}^{-1} \sum_{c=2}^{m} \binom{n_i}{2m-c} =  M \frac{(2m)!}{n_i(2m)} \sum_{c=2}^{m} \frac{n_i(2m-c)}{(2m-c)!} \leq M \sum_{c=2}^{m} \frac{n_i(2m-c)}{n_i(2m)} \\
    & \leq \frac{M}{(n_i-2m+2) (n_i-2m+1)} \leq \frac{M}{(n_0-2m+2) (n_0-2m+1)},
\end{align*}
where the last equality follows from \eqref{eq:compsamp}. Hence 
\[
n_0^2\left|\Xi^3_{i}\right| \leq M \frac{n_0^2}{(n_0-2m+2) (n_0-2m+1)}, 
\]
since \( n_0^2/\{(n_0-2m+2) (n_0-2m+1)\} \) is a strictly decreasing function of \( n_0 \) in its domain $n_0 \geq 2m$, its maximum is attained at  \( n_0=2m \), where it takes the value $2m^2$. Thus, since we are considering finite values of $m$, the result follows.
\end{proof}

Before proving Lemma \ref{lemma:VRk0} we introduce some useful bounds for $\V(L_{ik})$.

\begin{lemma}
Assume that \eqref{eq:independent_sample},\eqref{ass:boundenessh2},  \eqref{eq:compsamp} and \eqref{ass:boundsigma} hold. Then, for \( \forall i \), the following bounds hold:
\begin{enumerate}
    \item Under \( H_0 \),
    \begin{equation}
        M_1 / n_0^2 \leq \V(L_{ik}) \leq M_2 / n_0^2.
        \label{eq:lemma10.1gen}
    \end{equation}
    \item Under alternatives, 
    \begin{equation}
        V_k \geq M D_k/n_0. 
        \label{eq:lemma10LowerVk}
    \end{equation}
    \item If Assumption \ref{assumption}.i holds, for \( k > k_0 \), 
    \begin{equation}
     V_k \geq M_1 / n_0   , \quad V(L_{ik}) \leq M_2 / n_0.
    \label{eq:lemma10.2gen}
    \end{equation}
    \item If Assumption \ref{assumption}.ii holds and \( a_k / n_0 \leq M \), for large  $k$, 
    \begin{equation}
    a_k V_k \geq M_1 / n_0   , \quad a_k V(L_{ik}) \leq M_2 / n_0.
    \label{eq:lemma10.3gen}
    \end{equation}
    \item If Assumption \ref{assumption}.ii holds and \( a_k / n_0 \to \infty \), for large $k$, 
    \begin{equation}
     V_k \geq M_1 / n^2_0   , \quad V(L_{ik}) \leq M_2 / n^2_0.
    \label{eq:lemma10.4gen}
    \end{equation}

\end{enumerate}
\label{lemma:boundsVgen}
\end{lemma}

\begin{proof}

    \begin{enumerate}
    \item Suppose that \( H_0 \) holds. Then, by Remark \ref{remark}, as $\mathbb{L}_{ik2}=0$, we have that \( L_{ik} = \mathbb{L}_{i1} \). The lower bound follows from \eqref{ass:boundsigma} combined with \eqref{eq:lowerLi1}, while the upper bound is obtained from  \eqref{ass:boundenessh2} together with \eqref{eq:upperLi1}.  
    \item Note that $L_{ik}$ is a degree-$2m$ $U$-statistic with kernel
\[K(x_1,\ldots,x_{2m})=H(x_1,\ldots,x_{2m})-\theta_i^2-2\overline{\theta}\left\{\binom{2m}{m}^{-1}\sum_{1\leq j_1<\cdots<j_{m}\leq 2m}h(x_{j_1},\ldots,x_{j_m})-\theta_i\right\}.\]
Then, routine calculations show that $K_1(x) 
    =\left(\theta_i-\overline{\theta}\right)h_1(x)-\theta_i^2-\theta_i\overline{\theta}$, and hence 
$\V\{K_1(x)\}=(\theta_i-\overline{\theta})^2\sigma_1^2.$
This result, recalling \eqref{eq:varsigma}, implies that
$n_0 \V(L_{ik}) \geq M (\theta_i-\overline{\theta})^2\sigma_1^2$, and then, from \eqref{ass:boundsigma}, we get the result.
    
    \item In order to get the upper bound note that \( \E(L_{ik}) = \E(\mathbb{L}_{i1}) = \E(\mathbb{L}_{ik2}) = 0 \), by applying  \( c_r \) inequality, we get $n_0 \V(L_{ik}) \leq 2n_0 \V(\mathbb{L}_{i1}) + 2n_0 \V(\mathbb{L}_{ik2})$, and then, by recalling \eqref{eq:upperLi1} and \eqref{eq:upperLik2}, we obtain
$n_0 \V(L_{ik}) \leq M \E\{h^2(X_{i1},\ldots,X_{im})\}$.
From \eqref{ass:boundenessh2}, we conclude the upper bound.  
The lower bound follows from \eqref{eq:lemma10LowerVk} along with Assumption \ref{assumption}.i.

\item To establish the upper bound, we apply the \( c_r \) inequality, which gives  
$a_k n_0 \V(L_{ik}) \leq 2a_k n_0 \V(\mathbb{L}_{i1}) + 2a_k n_0 \V(\mathbb{L}_{ik2}).$ Using bound \eqref{eq:upperLi1} and noting that \( a_k/n_0 \) is bounded, we obtain  
$a_k n_0 \V(\mathbb{L}_{i1}) \leq M n_0^2 \V(\mathbb{L}_{i1}) \leq M.$ Similarly, from equation \eqref{eq:upperLik2theta} and the fact that \( a_k(\theta - \overline{\theta})^2 \) is uniformly bounded, we deduce that 
\[
a_k n_0 \V(\mathbb{L}_{ik2}) \leq M a_k (\theta_i - \overline{\theta})^2 \E\{h^2(X_{i1}, \ldots, X_{i2})\} \leq M \E\{h^2(X_{i1}, \ldots, X_{i2})\},
\]  
which is bounded from \eqref{ass:boundenessh2}.
Combining these results, we conclude the upper bound on \( a_k n_0\V(L_{ik}) \). To establish the lower bound, we consider \eqref{eq:lemma10LowerVk} along with assumption \ref{assumption}.ii. 

At this point myltiplying the previous inequality and recalling \eqref{eq:controlDk}, we get  
    \[
    a_kn_0V_k\geq M\left(1-\frac{\sqrt{a_k}}{n_0}\right) \geq M\left(1-\frac{1}{\sqrt{a_k}}\right),
    \]
    where the last inequality follows from $a_k/n_0 \leq M$. Since $a_k \to \infty$, the result holds for large $k$.  
    
 \item To get the lower bound, from Lemma \ref{lemma:variancegeneral} we have
\[V_k\geq M\left(\frac{D_k}{ n_0}-\frac{1}{k}\sum_{i=1}^k \left\{ \left|\theta_i-\overline{\theta}\right|\left|\Xi^3_{i}\right|+\V(\mathbb{L}_{i1})\right\}\right),\]
 
 Now we apply Lemma \ref{lemma:Xibounded} and \eqref{eq:lowerLi1}, which results in
\[V_k\geq M\left(\frac{D_k}{ n_0}-\frac{1}{n_0^2}\frac{1}{k}\sum_{i=1}^k \left|\theta_i-\overline{\theta}\right|+\frac{1}{n_0^2}\right),\]
and by applying moment monotonicity we get
\[
V_k \geq M\left(\frac{D_k}{ n_0} - \frac{D_k^{1/2}}{n_0^2} + \frac{1}{n_0^2}\right).
\]
Now we multiply the previous equation by $n_0^2$ and  note that \( n_0 D_k = (n_0/a_k) a_k D_k \to 0 \) since \( a_k D_k \) is bounded for $k \geq k_0$, \( n_0/a_k \to 0 \). Additionally, as \( D_k \to 0 \), the term \( D_k^{1/2} \) tends to zero. Then for sufficiently large $k$, we get the lower bound.

On the other hand, by applying $c_r$ inequality, it follows that $n_0^2 \V(L_{ik}) \leq  2n^2_0 \V(\mathbb{L}_{i1}) + 2n^2_0 \V(\mathbb{L}_{ik2}).$
From \eqref{ass:boundenessh2} and \eqref{eq:upperLi1}, $n^2_0 \V(\mathbb{L}_{i1})$ is bounded. From \eqref{eq:upperLik2theta}, we obtain
\[ n^2_0 \V(\mathbb{L}_{ik2}) \leq  M n_0(\theta_i-\overline{\theta})^2 \E\{h^2(X_{i1},\ldots, X_{im})\}=\frac{n_0}{a_k} a_k (\theta_i-\overline{\theta})^2 \E\{h^2(X_{i1},\ldots, X_{im})\},\]
which is bounded for large $k$ recalling that $n_0/a_k \to 0$, $a_k(\theta-\overline{\theta})^2$ is uniformly bounded and \eqref{ass:boundenessh2}. Then it follows that $n_0^2 \V(L_{ik})$ remains bounded.
    \end{enumerate}
\end{proof}

\begin{proof}[Proof of Lemma \ref{lemma:VRk0}]
    Note that,
\begin{equation}
   R_k = \frac{T_{k,\text{Lin}}}{k}  -R_{k1} + 2R_{k2}
\label{eq:decompRk} 
\end{equation}
where
\[
R_{k1} = \frac{1}{k^2} \sum_{i \neq j} (\widehat{\theta}_i - \theta_i)(\widehat{\theta}_j - \theta_j) \quad \text{and} \quad R_{k2} = \frac{1}{k^2} \sum_{i=1}^{k} (\theta_i - \overline{\theta})(\widehat{\theta}_i - \theta_i).
\]
So, to show the result it suffices to see that
\begin{equation}
    \V(R_{ki})/\V(T_{k, Lin}) \to 0,
    \label{eq:vrki}
\end{equation}
for $i=1,2$. We now upper-bound $\V(R_{k1})$ as follows:
\begin{align}
\V(R_{k1}) &= \E(R_{k1}^2) = \frac{1}{k^4} \sum_{i \neq j, a \neq b} \E \left\{ (\widehat{\theta}_i - \theta_i)(\widehat{\theta}_j - \theta_j)(\widehat{\theta}_a - \theta_a)(\widehat{\theta}_b - \theta_b) \right\} \notag \\
&= \frac{2}{k^4} \sum_{i \neq j} \V(\widehat{\theta}_i)\V(\widehat{\theta}_j) \leq \frac{2}{k^2} \left\{ \frac{1}{k} \sum_{i=1}^{k} \V(\widehat{\theta}_i) \right\}^2, 
\label{eq:VRk1}
\end{align}
where the third equality follows from the independence of samples. Specifically, the expectation 
\[
\E \left\{ (\widehat{\theta}_i - \theta_i)(\widehat{\theta}_j - \theta_j)(\widehat{\theta}_a - \theta_a)(\widehat{\theta}_b - \theta_b) \right\}=
\begin{cases}
    0, & \text{if } \{i,j\} \neq \{a,b\}, \\
    \V(\widehat{\theta}_i)\V(\widehat{\theta}_j), & \text{if } \{i,j\} = \{a,b\},
\end{cases}
\]
retains only terms where pairs of indices match, a situation that occurs twice, resulting in the expression \eqref{eq:VRk1}.

At this stage, by appliying \eqref{eq:varsigma} and recalling  \eqref{eq:sigmaorder} along with \(\sigma_{im}^2 \leq \E\left\{h^2(X_{i1},\ldots,X_{im})\right\} \), we can obtain an upper bound on \( \V(\widehat{\theta}_i) \), given by
\begin{align}
    \V(\widehat{\theta}_i) \leq \binom{n_i}{m}^{-1} \sum_{c=1}^{m} \binom{m}{c} \binom{n_i-m}{m-c} \sigma_{im}^2&\leq M \binom{n_i}{m}^{-1} \left\{\sum_{c=0}^{m} \binom{m}{c} \binom{n_i-m}{m-c}-\binom{n_i-m}{m}\right\} \nonumber \\
    &\leq  M  \left\{1-\binom{n_i}{m}^{-1}\binom{n_i-m}{m}\right\},
    \label{eq:boundvwtheta}
\end{align}
where the last inequality follows from Vandermonde's identity \eqref{eq:vandermonde}. In addition, we have that
\begin{align*}
    1-\binom{n_i}{m}^{-1}\binom{n_i-m}{m}&=\frac{n_i(n_i-1)\ldots(n_i-m+1)-(n_i-m)(n_i-m-1)\ldots(n_i-2m+1)}{n_i(n_i-1)\ldots(n_i-m+1)} \\
    &\leq \frac{\prod_{j=0}^{m-1} (n_i-j)-\prod_{j=0}^{m-1} (n_i-m-j)}{n_i^m}\leq \frac{n_i^m-(n_i-2m)^m}{n_i^m} \\
    &=n_i^{-m}\sum_{r=1}^{m}\binom{m}{r}(-1)^r n_i^{m-r} (2m)^r\leq (2m)^m \sum_{r=1}^m \binom{m}{r} n_i^{-r}\leq 2^{2m} m^{m+1}/n_i \leq M/n_0, 
\end{align*}
where the last inequality follows from \eqref{eq:compsamp}. Substituting the previous equation into \eqref{eq:boundvwtheta}, we have that
\begin{equation}
    V(\widehat{\theta}_i)\leq M/n_0,\;\text{ for } \forall i,
    \label{eq:boundvthetan0}
\end{equation}
and thus, from \eqref{eq:VRk1} we get the bound
\begin{equation}
    \V(R_{k1}) \leq M / (n^2_0 k^2).
    \label{eq:VRk1bound}
\end{equation}

Next, to study $R_{k1}$, we consider different cases.

If \( H_0 \) is true, from \eqref{eq:vtklin}, the lower bound in \eqref{eq:lemma10.1gen} and \eqref{eq:VRk1bound}, we get that 
\begin{equation}
    \V(R_{k1}) / \V(T_{k,\text{Lin}}) \leq M_1/k,
    \label{eq:VRk1/VTk}
\end{equation}
and thus $\V(R_{k1}) / \V(T_{k,\text{Lin}}) \to 0$.  

If $H_0$ does not hold, using \eqref{eq:vtklin} and \eqref{eq:lemma10LowerVk}, we get
\begin{equation}
    \V(T_{k,\text{Lin}}) \geq M D_k/(kn_0).
    \label{eq:Vtklinbound}
\end{equation}

\begin{itemize}
    \item On the one hand, if Assumption \ref{assumption}.i holds, then from \eqref{eq:lemma10.2gen}, $\V(T_{k,\text{Lin}})\geq M_1/(n_0 k),$ for $k > k_0$. This fact along with \eqref{eq:VRk1bound} gives that \eqref{eq:vrki} is true for $i=1$. 
    \item On the other hand, if Assumption \ref{assumption}.ii holds, then from \eqref{eq:lemma10.3gen} $a_k \V(T_{k,\text{Lin}}) \geq M_1/(n_0 k),$ for  $k > k_0$ when  $a_k/n_0$ remains bounded, and the same applies when $a_k/n_0 \to \infty$ since, from \eqref{eq:lemma10.4gen}
    $a_kV_k \geq M_1 a_k/n^2_0 \geq M_1/n_0$  for large $k$. Hence using \eqref{eq:VRk1bound}, we get 
    \begin{equation}
      a_k \V(R_{k1})/\{a_k \V(T_{k,\text{Lin}})\} \leq k^{\eta} \V(R_{k1})/\{a_k \V(T_{k,\text{Lin}})\} \leq M_1/(n_0 k^{1-\eta}) \to 0.
      \label{eq:Caso3Rk1}
    \end{equation}
    
\end{itemize}

For $R_{k2}$, we have that
\[\V(R_{k2})=\frac{1}{k^4} \sum_{i=1}^k (\theta_i-\overline{\theta})^2 \V(\widehat{\theta_i}).\]
So introducing \eqref{eq:boundvthetan0} in the previous equation, we get
\begin{equation}
    \V(R_{k2})\leq M D_k/(n_0k^3).
    \label{eq:Vrk2bound}
\end{equation}
If $H_0$ is true, then $\V(R_{k_2})$ is trivially 0 and \eqref{eq:vrki} holds for $i=2$. Under alternatives, from \eqref{eq:Vtklinbound} and \eqref{eq:Vrk2bound}, we get 
\[\V(R_{k2})/\V(T_{k, Lin})\leq M/k^2, \]
and hence, \eqref{eq:vrki} holds for $i=2$.
\end{proof}

\begin{proof}[Proof of Lemma \ref{lemma:Tlinlimit}]
To prove the result, we verify Lindeberg's condition
\begin{equation}
    \ell(\varepsilon) = \frac{\frac{1}{k} \sum_{i=1}^k \E \left[ L_{ik}^2 \, I\{ L^2_{i,k} > \varepsilon \sum_{i=1}^k \V(L_{ik}) \} \right]}{\frac{1}{k} \sum_{i=1}^k \V(L_{ik})} \rightarrow 0, \quad \forall \, \varepsilon > 0,
    \label{eq:lindeberg}
\end{equation}
where \( I(\cdot) \) denotes the indicator function.

We proceed by analyzing different cases, starting with the assumption that \( H_0 \) holds. First, we multiply and divide \(\ell(\varepsilon)\) by \(n_0^2\). Then, invoking the lower bound in \eqref{eq:lemma10.1gen}, we obtain that, to prove \eqref{eq:lindeberg}, it suffices to show that  
\begin{equation}
\mu_i(\varepsilon) = \E \left[ n_0^2 L_{ik}^2 \, I \left\{ \, n_0^2 L_{ik}^2 > \varepsilon \sum_{i=1}^k \V(n_0 L_{ik}) \right\} \right] \rightarrow 0, \quad \forall \, \varepsilon > 0
    \label{eq:lindebergn2}
\end{equation}
Next, observe that
\[
I\Bigg\{ n_0^2 L_{ik}^2 > \varepsilon \sum_{i=1}^k \V(n_0 L_{ik}) \Bigg\} \leq I( n_0^2 L_{ik}^2 > \varepsilon k M_1 ).
\]
Consequently, we obtain the bound
\begin{equation}
    \mu_i(\varepsilon) \leq \E \left\{ n_0^2 L_{ik}^2 \, I( n_0^2 L_{ik}^2 > \varepsilon k M_1 ) \right\} \leq \E(n_0^2 L_{ik}^2).
    \label{eq:mui1bound}
\end{equation}
Since \( \E( L_{ik}) = 0 \), from the upper bound in \eqref{eq:lemma10.1gen} it follows that $\E(n_0^2 L_{ik}^2)$ is finite. Using this result and applying the Dominated Convergence Theorem to \eqref{eq:mui1bound}, we deduce that for each \( i \),
\[
\lim\mu_i(\varepsilon) \leq \lim \E \left\{ n_0^2 L_{ik}^2 \, I\left( n_0^2 L_{ik}^2 > \varepsilon k M_1 \right) \right\} 
= \E \left\{ n_0^2 L_{ik}^2 \, \lim I\left( n_0^2 L_{ik}^2 > \varepsilon k M_1 \right) \right\} = 0, \quad \forall \varepsilon>0.
\]
Therefore, \eqref{eq:lindeberg} holds when \( H_0 \) is true.

Now assume that Assumption~\ref{assumption}.i holds. Recalling the lower bound in \eqref{eq:lemma10.2gen}, it follows that, to establish~\eqref{eq:lindeberg}, it suffices to verify that
\[
\mu_i(\varepsilon) = \E \left[ n_0 L^2_{ik} \, I\left\{ n_0 L^2_{ik} > \varepsilon \sum_{i=1}^k \V\left( \sqrt{n_0} L_{ik} \right) \right\} \right] \longrightarrow 0, \quad \text{for all } \varepsilon > 0.
\]
From the lower bound in \eqref{eq:lemma10.2gen}, it follows that
\[I\Bigg\{ n_0 L^2_{ik} > \varepsilon \sum_{i=1}^k \V(\sqrt{n_0} L_{ik})\Bigg\} \leq I( n_0 L_{ik}^2 > \varepsilon k M_1 ),\] 
which yields
\begin{equation}
    \mu_i(\varepsilon) \leq \E \left\{ n_0 L_{ik}^2 \, I( n_0 L_{ik}^2 > \varepsilon k M_1 ) \right\} \leq \E(n_0 L_{ik}^2).
    \label{eq:muiass1ibound}
\end{equation}
Once again, since \( \E(L_{ik}) = 0 \), the upper bound in \eqref{eq:lemma10.2gen} implies that \( \E(n_0 L_{ik}^2) \) is bounded. This allows us to apply the Dominated Convergence Theorem to \eqref{eq:muiass1ibound}, which in turn confirms that \eqref{eq:lindeberg} holds under Assumption \ref{assumption}.i.  

Suppose now that Assumption \ref{assumption}.ii holds and that $a_k/n_0 \leq M$. From the lower bound in  \eqref{eq:lemma10.3gen}, it suffices to show that 
\[\mu_i(\varepsilon)=\E \left[  n_0 a_k L^2_{ik} \, I\Bigg\{ a_k n_0 L^2_{ik} > a_k \varepsilon \sum_{i=1}^k \V(\sqrt{n_0} L_{ik})\Bigg\} \right] \rightarrow 0, \quad \forall \, \varepsilon > 0.\]
In this case, from the lower bound in \eqref{eq:lemma10.3gen}, we get   
\[I\Bigg\{ a_k n_0 L^2_{ik} > a_k \varepsilon \sum_{i=1}^k \V(\sqrt{n_0} L_{ik})\Bigg\} \leq I( a_k n_0 L_{ik}^2 > \varepsilon k M_1),\] 
and then
\begin{equation}
    \mu_i(\varepsilon) \leq \E \left\{ a_k n_0 L_{ik}^2 \, I\left( a_k n_0 L_{ik}^2 > \varepsilon k M_1 \right) \right\} \leq \E\left(a_k n_0 L_{ik}^2\right).
    \label{eq:muiass12bound}
\end{equation}
As $\E(L_{ik}) = 0$, from the upper bound in \eqref{eq:lemma10.3gen}, $\E\left(a_k n_0 L_{ik}^2\right)<M$.
Hence, we apply the Dominated Convergence Theorem noting that since $a_k$ is positive and $a_k\leq k^\eta$, for some $0<\eta<1$, we have
\[\lim  I\left( a_k n_0 L_{ik}^2 > \varepsilon k M_1 \right)=\lim  I\left( n_0 L_{ik}^2 > \varepsilon k/a_k M_1 \right) \leq \lim  I\left( n_0 L_{ik}^2 > \varepsilon k^{1-\eta} M_1 \right)=0.
\]

Finally, we assume that Assumption \ref{assumption}.ii holds and that \( a_k / n_0 \to \infty \). Recalling the lower bound in \eqref{eq:lemma10.4gen}, we conclude that, to establish \eqref{eq:lindeberg}, it suffices to verify that \eqref{eq:lindebergn2} holds.  Now, from the upper bound in \eqref{eq:lemma10.4gen}, we get that $\E(n_0^2L_{ik}^2)$ is finite.
Following the same reasoning as in the first case, one gets that \eqref{eq:lindeberg} holds.  
\end{proof}

\begin{proof}[Proof of Theorem \ref{Theorem:Tklim2}]
Using Lemma \ref{lemma:VRk0} and Chebyshev's inequality, for all $\varepsilon > 0$, we have
\[
\Pr\Big(\Big|\frac{R_k}{\sqrt{\V(T_{k, \text{Lin}})}}\Big| > \varepsilon \Big) \leq \frac{\V(R_k)}{\varepsilon^2 \V(T_{k, \text{Lin}})} \to 0,
\]
which implies that $ R_k /\sqrt{\V(T_{k, \text{Lin}})} \tp 0$. With this result, recalling \eqref{eq:decompTk}, Lemma \ref{lemma:Tlinlimit}, and using Slutsky's Theorem, it follows that
\[
\frac{T_k - D_k}{\sqrt{\V(T_{k, \text{Lin}})}} \tl Z,
\]
where $Z$ follows a standard normal distribution, completing the proof.
\end{proof}

\begin{proof}[Proof of Proposition \ref{proposition:varianceestimator}]
In first place, $S^2_{\hat{l}}$ can be decomposed as \begin{equation}
S^2_{\hat{l}}=S_1+S_2+S_3+S_{12}+S_{13}+S_{23},
\label{eq:decompS}
\end{equation} where
\[S_1=\frac{1}{k-1}\sum_{i=1}^k (\hat{l}_{ik}-l_{ik})^2, \quad S_2=\frac{k}{k-1} (\overline{l}-\widehat{\overline{l}})^2, \quad S_3=\frac{1}{k-1}\sum_{i=1}^k (l_{ik}-\overline{l})^2, \]
and $S_{ij}^2 \leq S_iS_j$ for $i,j=1,2,3$. By definition, $S_1$ and $S_2$ can be written as follows
\begin{equation}
    S_1=\frac{4}{k-1} \sum_{i=1}^k \widehat{\theta_i}^2\Bigg\{\frac{1}{k}\sum_{j=1}^k(\widehat{\theta}_i-\theta_i)\Bigg\}^2, \quad S_2=\frac{4k}{k-1}\Bigg(\frac{1}{k}\sum_{i=1}^k\widehat{\theta_i}\Bigg)^2 \Bigg\{\frac{1}{k}\sum_{j=1}^k(\widehat{\theta}_i-\theta_i)\Bigg\}^2
    \label{eq:S1S2}
\end{equation}
As we are assuming \eqref{ass:boundenessh2}, we have, as we saw in the proof of Lemma \ref{lemma:consistency}, that $\E\left(\widehat{\theta_i}^2\right)<M$ and hence  \begin{equation}
    \frac{1}{k-1}\sum_{i=1}^k\widehat{\theta_i}^2=O_{\Prob}(1).
    \label{eq:k-1widtheta2}
\end{equation}
Similarly, as $\V(\widehat{\theta_i})$ is finite under \eqref{ass:boundenessh2}, it follows that
\begin{equation}
    \frac{1}{k}\sum_{i=1}^k\widehat{\theta_i}=O(1).
    \label{eq:k-1widtheta}
\end{equation} Now we check Lindeberg's condition 
\begin{equation}
    \ell(\varepsilon) = \frac{\frac{1}{k} \sum_{i=1}^k \E \left[ (\widehat{\theta}_i-\theta_i)^2 \, I\Bigg\{ (\widehat{\theta}_i-\theta_i)^2 > \varepsilon \sum_{i=1}^k \V(\widehat{\theta}_i) \Bigg\} \right]}{\frac{1}{k} \sum_{i=1}^k \V(\widehat{\theta}_i)} \rightarrow 0, \quad \forall \, \varepsilon > 0,
    \label{eq:lindebergthetai}
\end{equation}
Now, we observe that, from \eqref{ass:boundsigma}, \( \sigma^2_{i1} \) is bounded  away from zero, and thus \( \widehat{\theta_i} \) is a non-degenerate \( U \)-statistic. Therefore, from \eqref{eq:varsigma} and \eqref{eq:compsamp}, we conclude that 
\begin{equation}
n_0 \V(\widehat{\theta}_i) > M \quad \text{ for } i = 1, \ldots, k,
\label{eq:boundVthetaresultadovarianza}
\end{equation}
and thus \( (1/k) \sum_{i=1}^k n_0 \V(\widehat{\theta}_i) \) is bounded from below by some constant.
Thus, to check Lindeberg's condition it suffices to see that
\[\mu_i(\varepsilon)=\E \left[  n_0  (\widehat{\theta}_i-\theta_i)^2 \, I\{  n_0 (\widehat{\theta}_i-\theta_i)^2 >  \varepsilon \sum_{i=1}^k n_0\V(\widehat{\theta_i})\} \right] \rightarrow 0, \quad \forall \, \varepsilon > 0.\]
Hence, considering the bound \eqref{eq:boundVthetaresultadovarianza}, it follows that
\begin{equation*}
    I\left\{ n_0 (\widehat{\theta}_i-\theta_i)^2 >  \varepsilon \sum_{i=1}^k n_0\V(\widehat{\theta_i}) \right\} 
    \leq I\left\{ n_0 (\widehat{\theta}_i-\theta_i)^2 >  \varepsilon M k \right\},
\end{equation*}
which implies
\begin{equation*}
    \mu_i(\varepsilon) \leq \E \left[ n_0  (\widehat{\theta}_i-\theta_i)^2 \, I\Bigg\{ n_0 (\widehat{\theta}_i-\theta_i)^2 >  \varepsilon M k \Bigg\} \right].
\end{equation*}
Now, as we showed in \eqref{eq:boundvthetan0}, as since \eqref{ass:boundenessh2} and \eqref{eq:compsamp} hold,  we have that 
\begin{equation*}
    \E\{ n_0  (\widehat{\theta}_i-\theta_i)^2\}=n_0 \V(\widehat{\theta}_i) \leq M.
\end{equation*}
Now, by applying the Dominated Convergence Theorem as in the proof of Lemma \ref{lemma:Tlinlimit}, we conclude that Lindeberg's condition is satisfied. This implies that
\begin{equation}
    \sqrt{n_0 k}\frac{ (1/k)\sum_{i=1}^k (\widehat{\theta}_i-\theta_i)}{\sqrt{(1/k)\sum_{i=1}^k n_0 \V(\widehat{\theta}_i)}} \tl Z,
    \label{eq:limitwidtheta}
\end{equation}
where $Z$ follows a standard normal distribution. Finally, since
    $M_1 \leq (1/k)\sum_{i=1}^k n_0 \V(\widehat{\theta}_i) \leq M_2$
and from \eqref{eq:limitwidtheta}, we obtain that 
\begin{equation}
    \frac{1}{k}\sum_{i=1}^k(\widehat{\theta}_i-\theta_i) = O_{\Prob}\left(\frac{1}{\sqrt{n_0 k}}\right).
    \label{eq:difwidthetatheta}
\end{equation}
Finally, from \eqref{eq:k-1widtheta2}, \eqref{eq:k-1widtheta} and \eqref{eq:difwidthetatheta}, it follows that 
\begin{equation*}
    S_i = O_{\Prob}\left(\frac{1}{n_0 k}\right), \quad i=1,2.
    \label{eq:Si0}
\end{equation*}

Next, we are going to prove that $S_i/\E(S^2_l)=o_{\Prob}(1), \; i=1,2$. 

First, it is important to note that  \eqref{ass:boundenessh2} is sufficient to guarantee that \( \V(L_{ik}) \) is upper-bounded, as shown in Lemma \ref{lemma:boundsVgen}. Consequently, \( V_k \) is also upper-bounded.
\begin{itemize}
\item If \( H_0 \) holds, from the lower bound in \eqref{eq:lemma10.1gen}, since in this case \( \Delta_k = 0 \), we have
\[
\frac{S_i}{\E(S^2_l)} = \frac{S_i}{V_k}=O_{\Prob}(n_0/k) = o_{\Prob}(1),
\]  
for $i=1,2$, because it is assumed that \( n_0/k \to 0 \).
\item Consider that Assumption \ref{assumption}.i holds, recalling the lower bound in in \eqref{eq:lemma10.2gen}, we have that $ V_k=O(1/n_0)$, therefore \[\frac{S_i}{\E(S^2_l)}\leq \frac{S_i}{V_k} = O_{\Prob}(1/k) = o_{\Prob}(1),\]
for $i=1,2$.
\item Under Assumption \ref{assumption}.ii and $a_k/n_0 \leq M$, recalling the lower bound in in \eqref{eq:lemma10.3gen}, it holds that
\[
\frac{a_k S_i}{a_k\E(S^2_l)}\leq \frac{a_k S_i}{a_k V_k} = O_{\Prob}(1/k^{1-\eta}) = o_{\Prob}(1),
\] 
for $i=1,2$.
\item Under Assumption \ref{assumption}.ii and the condition \(a_k/n_0 \to \infty\), recalling the lower bound in \eqref{eq:lemma10.4gen}, we get  
\[
\frac{ S_i}{\E(S^2_l)} \leq \frac{S_i}{ V_k} = O_{\Prob}(n_0/k) = o_{\Prob}(1),
\] 
for \(i=1,2\), following the same reasoning as in the \(H_0\) case.
\end{itemize}

On the other hand, $L_{ik}=l_{ik}-\E(l_{ik})$ and $T_{k,\text{Lin}}=\overline{l}-\E(\overline{l})$. Hence, $S_3$ can be expressed as
\begin{align*}
    S_3 &= \frac{1}{k-1}\sum_{i=1}^k\Big\{L_{ik}-T_{k,\text{Lin}}+\E(l_{ik})-\E(\overline{l})\Big\}^2 \\
    &= \frac{1}{k-1}\sum_{i=1}^k L_{ik}^2-\frac{k}{k-1} T^2_{k,\text{Lin}} +\underbrace{\frac{1}{k-1}\sum_{i=1}^k \{\E(l_{ik})-\E(\overline{l})\}^2}_{=\Delta_k} + \frac{2}{k-1} \sum_{i=1}^k \{\E(l_{ik})-\E(\overline{l})\}L_{ik}.
\end{align*}
Lemma \ref{lemma:Tlinlimit} establishes that, under either $H_0$ or Assumption \ref{assumption},   
\[
\sqrt{k} \frac{T_{k,\text{Lin}}}{\sqrt{V_k}} \tl Z.
\]
As a consequence, it follows that  
\[
\frac{T_{k,\text{Lin}}}{\sqrt{V_k}} = O_{\Prob}\left(\frac{1}{\sqrt{k}}\right),
\]
which further implies that  
\[
\left|\frac{T_{k,\text{Lin}}}{\sqrt{\E(S^2_l)}}\right| \leq \frac{|T_{k,\text{Lin}}|}{\sqrt{V_k}} = O_{\Prob}\left(\frac{1}{\sqrt{k}}\right).
\]
Finally, this leads to  
\[
\frac{T^2_{k,\text{Lin}}}{\E(S^2_l)} = O_{\Prob}\left(\frac{1}{k}\right) = o_{\Prob}(1).
\]
We now compute
\begin{align*}
&\E\left(\left[\frac{1}{k-1} \frac{\sum_{i=1}^k \{\E(l_{ik})-\E(\overline{l})\}L_{ik}}{\E(S^2_l)}\right]^2\right) \notag \\
&\quad = \frac{1}{(k-1)^2} \frac{\sum_{i=1}^k \{\E(l_{ik})-\E(\overline{l})\}^2 \E(L^2_{ik})
+ \sum_{r\neq s} \{\E(l_{rk})-\E(\overline{l})\} \{\E(l_{sk})-\E(\overline{l})\} \E(L_{rk} L_{sk})}{\E^2(S^2_l)} \notag \\
& \quad = \frac{1}{(k-1)^2} \sum_{i=1}^k \frac{\{\E(l_{ik})-\E(\overline{l})\}^2}{\E(S^2_l)} \frac{\E(L^2_{ik})}{\E(S^2_l)}\leq \frac{1}{(k-1)^2} \sum_{i=1}^k \frac{\{\E(l_{ik})-\E(\overline{l})\}^2}{\E(S^2_l)} \frac{\V(L_{ik})}{V_k},
\end{align*}
where we used that \(\E(L_{rk} L_{sk}) = \E(L_{rk})\E(L_{sk})\) due to the independence of the random samples, along with the fact that \(\E(L_{ik}) = 0, \; \forall i\). Now, from the bounds in Lemma \ref{lemma:boundsVgen}, we get
\begin{equation}
    \E\left(\left[\frac{1}{k-1} \frac{\sum_{i=1}^k \{\E(l_{ik})-\E(\overline{l})\}L_{ik}}{\E(S^2_l)}\right]^2\right)\leq \frac{M}{k-1}\frac{\Delta_k}{\E[S^2_l]}\leq \frac{M}{k-1},
    \label{eq:esperanzadeltalik}
\end{equation}
since $\Delta_k\leq \E[S^2_l]$.
Therefore, by applying Markov’s inequality, we obtain  
\[
\frac{1}{k-1} \frac{\sum_{i=1}^k \{\E(l_{ik}) - \E(\overline{l})\}L_{ik}}{\E(S^2_l)} = O_{\Prob}\left(\frac{1}{\sqrt{k}}\right) = o_{\Prob}(1).
\]
We have proven that
\begin{equation}
    \frac{S_3 - \E(S^2_l)}{\E(S^2_l)}= \frac{\frac{1}{k-1} \sum_{i=1}^k L^2_{ik} - V_k}{\E(S^2_l)} + o_{\Prob}(1).
    \label{eq:variancefinalconvergence}
\end{equation}
Under \( H_0 \), we apply Markov's inequality, the \( c_r \) inequality, and a refinement of the $c_r$ inequality from Theorem 2 of von Bahr (1965), obtaining
\begin{align}
    \Pr\left[\left| \frac{1}{k} \sum_{i=1}^k \left\{ n_0^2 L^2_{ik} - \V(n_0 L_{ik}) \right\} \right| > \varepsilon \right] 
    &\leq \frac{1}{\varepsilon^{1+\delta/2}} \frac{\E \left[ \left| \sum_{i=1}^k \left\{ n_0^2 L^2_{ik} - \V(n_0 L_{ik}) \right\} \right|^{1+\delta/2} \right]}{k^{1+\delta/2}} \nonumber \\
    &\leq \frac{2}{k^{\delta/2} \varepsilon^{1+\delta/2}} \frac{1}{k} \sum_{i=1}^k \E \left\{ \left| n_0^2 L^2_{ik} - \V(n_0 L_{ik}) \right|^{1+\delta/2} \right\} \nonumber \\
    &\leq \frac{4}{k^{\delta/2} \varepsilon^{1+\delta/2}} \frac{1}{k} \sum_{i=1}^k \left\{ \E \left( \left| n_0 L_{ik} \right|^{2+\delta} \right) + \V^{1+\delta/2}(n_0 L_{ik}) \right\} \nonumber \\
    &\leq  \frac{M}{k^{\delta/2} \varepsilon^{1+\delta/2}},
    \label{eq:convergenceL2Vk}
\end{align}
The last inequality follows from the upper bound in \eqref{eq:lemma10.1gen} and the assumption that \(\E \left( \left| n_0 L_{ik} \right|^{2+\delta} \right)\) is uniformly bounded.
This establishes that
\[
\frac{1}{k} \sum_{i=1}^k \left\{ n_0^2 L^2_{ik} - \V(n_0 L_{ik}) \right\} = o_{\Prob}(1),
\]
which, recalling the lower bound on  $V_k$  in \eqref{eq:lemma10.1gen}, implies
\[
\frac{\frac{1}{k} \sum_{i=1}^k L^2_{ik} - V_k}{V_k} = o_{\Prob}(1).
\]
Thus, from \eqref{eq:variancefinalconvergence}, we conclude that under \( H_0 \),
\begin{equation}
    \frac{S_3 - \E(S_l)}{\E(S^2_l)} = o_{\Prob}(1).
    \label{eq:S3ratioconsistent}
\end{equation}
If Assumption \ref{assumption}.i holds, recalling \eqref{eq:lemma10.2gen} and following a similar reasoning, we obtain that \eqref{eq:S3ratioconsistent} holds. Otherwise, if Assumption \ref{assumption}.ii holds and \( a_k/n_0 \) is bounded, applying the bound \eqref{eq:lemma10.3gen}, and noting that  
\[
\E\left[\left|\sqrt{a_k n_0}L^2_{ik}\right|^{2+\delta}\right] = \left|\frac{a_k}{n_0}\right|^{1+\delta/2} \E\left[\left| n_0L^2_{ik}\right|^{2+\delta}\right] \leq M_1 \E\left[\left| n_0L^2_{ik}\right|^{2+\delta}\right],
\]  
proceeding analogously, we conclude that \eqref{eq:S3ratioconsistent} also holds. 
In the case where Assumption \ref{assumption}.ii holds and $a_k/n_0 \to \infty$, from the bounds in \eqref{eq:lemma10.4gen}, and following the same steps as in the first case, \eqref{eq:S3ratioconsistent} follows.

With the above results, from \eqref{eq:S3ratioconsistent} and  \( S_i = o_{\Prob}(1) \) for \( i = 1,2 \), we conclude that  
\[
    \frac{S^2_l}{\E(S^2_l)} = 1 - \frac{S^2_l - \E(S_l)}{\E(S^2_l)} = 1 + o_{\Prob}(1),
\]  
which completes the proof.
\end{proof}

\begin{lemma}
    Suppose that \eqref{ass:boundenessh2} holds. Then, it follows that \[ \Delta_k \leq M \frac{k}{k-1} D_k,\]
    where $\Delta_k$ is defined in \eqref{eq:Deltak}.
    \label{lemma:deltabounded}
\end{lemma}
\begin{proof}
    The proof follows a similar argument to that of Lemma 7 in \cite{JimenezGamero2025}.
\end{proof}

\begin{proof}[Proof of Theorem \ref{theorem:power}]
    First, the critical region is given by
    \[\mathcal{T}_k > z_{1-\alpha} \Longleftrightarrow \sqrt{k}\frac{ T_k}{S_{\hat{l}}}> z_{1-\alpha} \Longleftrightarrow \sqrt{k}T_k > z_{1-\alpha}S_{\hat{l}}.\] 

From Theorem  \ref{Theorem:Tklim2}, under Assumption \ref{assumption}, we have that 
\[0 \leq \sqrt{k}T_k = \sqrt{k}D_k + \sqrt{V_k} Z_k, \quad Z_k \tl Z. \]
On the other hand, from Proposition \ref{proposition:varianceestimator}, we know that
\[0\leq S_{\hat{l}}=\sqrt{V_k+\Delta_k} (1+s_k),\]
 where $s_k=o_{\Prob}(1)$. With these results, the critical region $\mathcal{T}_k> z_{1-\alpha}$ can be written as
\begin{equation}
    \sqrt{k} D_k + \sqrt{V_k}Z_k > z_{1-\alpha}\sqrt{V_k + \Delta_k}(1+s_k) .
    \label{eq:criticalregion}
\end{equation}
Now, we analyze the behavior of both sides of the inequality case by case:

\begin{enumerate}[label=\alph*), leftmargin=*]  
\item If Assumption \ref{assumption}.i holds, for the left-hand side of \eqref{eq:criticalregion}, recalling the upper-bound of \eqref{eq:lemma10.2gen}, we have that
\[\sqrt{k} D_k + \sqrt{V_k}Z_k\geq M\sqrt{k}  - \sqrt{V_k}|Z_k|\geq M\left(\sqrt{k}-\frac{|Z_k|}{n_0}\right) \tp \infty, \]
for either bounded $n_0$ or $n_0 \to \infty$. To prove the result, it suffices to show that \( S_{\hat{l}} \) is upper-bounded, which reduces to establishing the same for \( \Delta_k \). From Lemma \ref{lemma:deltabounded}, it is sufficient to verify that \( D_k \) is upper-bounded, which is immediate from \eqref{eq:thetaacotado}. Then, all the above implies that \( P(\mathcal{T}_k > z_{1-\alpha}) \to 1 \).

\item
If Assumption \ref{assumption}.ii holds and $\sqrt{k}/a_k \to \infty$,  for the left-hand side of \eqref{eq:criticalregion}, we get
\[\sqrt{k}D_k+\sqrt{V_k} Z_k, \geq M \left(\frac{\sqrt{k}}{a_k}-|Z_k|\right) \tp \infty,\]
where we have used that  $V_k\leq M$ from Lemma \ref{lemma:boundsVgen}. To prove the result it suffices to check that  $\Delta_k$ is upper-bounded, which follows from the same argument as in the previous case.
\item If Assumption \ref{assumption}.ii holds, $M_1\leq \sqrt{k}/a_k \leq M_2$ and $n_0 \to \infty$, for the left-hand side of \eqref{eq:criticalregion} we have that
\[\sqrt{k}D_k+\sqrt{V_k} Z_K\geq M \left(\sqrt{\frac{k}{a_k}}-\frac{|Z_k|}{\sqrt{n_0}}\right) \tp M>0,\]
where we have used that from Lemma \ref{lemma:boundsVgen}, $V_k \leq M/n_0$ in this case. As $V_k \to 0$, to prove the result it suffices to check that $\Delta_k \to 0$, which holds from Lemma \ref{lemma:deltabounded} because $D_k \to 0$.

\item If Assumption~\ref{assumption}.ii holds, \( \sqrt{k}/a_k \to 0 \), \( n_0 \to \infty \), and \( a_k/n_0 \leq M\), we express the critical region \eqref{eq:criticalregion} as
\[\sqrt{a_kn_0}\left(D_k+\sqrt{\frac{V_k}{k}} Z_k\right) > \sqrt{a_kn_0}\left\{z_{1-\alpha}\sqrt{\frac{V_k + \Delta_k}{k}}(1+s_k)\right\} .\]
For the left-hand side of the previous inequality, we get
\[\sqrt{a_kn_0}D_k+\sqrt{\frac{a_kn_0V_k}{k}} Z_k 
\geq M\left(1-\frac{|R_k|}{\sqrt{k}}\right) \tp M >0,\]
where we have used the upper bound on $V_k$ obtained in \eqref{eq:lemma10.3gen}. Sach a bound implies that $a_k n_0 V_k/k \to 0$ Then, to prove the result, it suffices to check that 
$\sqrt{a_k n_0 \Delta_k / k} \to 0.$
This follows from Lemma \ref{lemma:deltabounded}, together with the bound $a_k D_k \leq M$ stated in \eqref{eq:controlDk}, and the assumption $n_0 / k \to 0$ given in Proposition \ref{proposition:varianceestimator}.

\item If Assumption \ref{assumption}.ii holds, \( \sqrt{k}/a_k \to 0 \), \( n_0 \to \infty \), \( a_k/n_0 \to \infty\) and $\sqrt{k}n_0/a_k \to \infty$, we express the critical region as
\[a_k\left(D_k+\sqrt{\frac{V_k}{k}} Z_k\right) > a_k\left\{z_{1-\alpha}\sqrt{\frac{V_k + \Delta_k}{k}}(1+s_k)\right\} .\]
For the left-hand side of the expression above we have that
\[a_kD_k+\frac{a_k}{\sqrt{k}}\sqrt{V_k}Z_k \geq M\left(1-\frac{a_k}{n_0\sqrt{k}}|Z_k| \right)\tp M >0 ,\]
where we have used the upper bound on $V_k$ deduced from \eqref{eq:lemma10.4gen}. To get the result it suffices to check that $a_k^2\Delta_k/k$ tends to zero, which follows from Lemma \ref{lemma:deltabounded}, \eqref{eq:controlDk} and the fact that $a_k \leq k^{\eta}$.

\end{enumerate}
\end{proof}

\begin{proof}[Proof of Proposition \ref{proposition:conditioncfixedn0}]
    Since \(z_{1-\alpha} = -z_{\alpha}\), the power of the test is given by, 
\begin{align*}
\Prob\left(\mathcal{T}_k > z_{1-\alpha}\right) 
&= \Prob\left(\sqrt{k} \frac{T_k - D_k}{S_{\hat{l}}} > -z_{\alpha} - \sqrt{k} \frac{D_k}{S_{\hat{l}}}\right) \\
&= \Prob\left( \sqrt{k} \frac{T_k - D_k}{\sqrt{V_k}}  > -z_{\alpha}\frac{S_{\hat{l}}}{\sqrt{V_k}} - \sqrt{k} \frac{D_k}{\sqrt{V_k}}\right) \\
&= \Prob\left( \sqrt{k} \frac{T_k - D_k}{\sqrt{V_k}} \frac{\sqrt{V_k+\Delta_k}}{S_{\hat{l}}}   > -z_{\alpha} \sqrt{1+\frac{\Delta_k}{V_k}} - \sqrt{k} \frac{D_k}{\sqrt{V_k}} \frac{\sqrt{V_k+\Delta_k}}{S_{\hat{l}}}\right),
\end{align*}

By Proposition \eqref{proposition:varianceestimator} we have that that $\sqrt{V_k+\Delta_k}/S_{\hat{l}}=1+o_{\Prob}(1)$, and by Theorem \ref{Theorem:Tklim2} we know that $\Prob\left(\sqrt{k} \frac{T_k - D_k}{\sqrt{V_k}}<x\right)=\Phi(x)+o(1)$. Therefore, we get
\begin{equation}
    \Prob\left(\mathcal{T}_k > z_{1-\alpha}\right)\approx \Phi\left(\sqrt{k} \frac{D_k}{\sqrt{V_k}} + z_{\alpha} \sqrt{1+\frac{\Delta_k}{V_k}} \right).
    \label{eq:asymptoticpower}
\end{equation}
\end{proof}
\begin{enumerate}
    \item First, note that since $n_0$ is fixed, by Lemma \ref{lemma:boundsVgen}, $M_1\leq V_k\leq M_2$. We also know, from Lemma \ref{lemma:deltabounded}, that $\Delta_k \to 0$. These arguments imply that $\Delta_k/V_k \to 0$, and hence 
\begin{equation}
    \Prob\left(\mathcal{T}_k > z_{1-\alpha}\right)\approx \Phi\left(\sqrt{k} \frac{D_k}{\sqrt{V_k}} + z_{\alpha} \right).
    \label{eq:powersimplified}
\end{equation}
Now, from \eqref{eq:controlDk} and the aforementioned bounds for $V_k$, we have that
\begin{equation}
    M_1\frac{\sqrt{k}}{a_k}\leq\sqrt{k} \frac{D_k}{\sqrt{V_k}} \leq M_1\frac{\sqrt{k}}{a_k}.
    \label{eq:eq:boundssignalnoise}
\end{equation}
As we are assuming that $M_1\leq \sqrt{k}/a_k \leq M_2$, from the monotonicity of the cumulative distribution function, it follows that $\alpha < \mathcal{P} > z_{1-\alpha}) < 1$. If additionally, \( \sqrt{k} D_k / \sqrt{V_k} \to d \in (0, \infty) \), from \eqref{eq:powersimplified} we obtain that $\mathcal{P} = \Phi(z_{\alpha} + d)$.


\item We proceed as in the first case. First, it is important to note that
\[\frac{n_0/\sqrt{k}}{\sqrt{k}n_0/a_k}= \frac{a_k}{k}\leq k^{\eta-1}.\]
Next, we study the term $\Delta_k/V_k$, recalling Lemma \ref{lemma:deltabounded}, \eqref{eq:controlDk}, \eqref{eq:lemma10.4gen} and the equation above, we have that
\begin{equation}
    \frac{\Delta_k}{V_k}\leq M\frac{k}{k-1}\frac{n_0^2}{a_k}= M\frac{k}{k-1}\frac{\sqrt{k}n_0}{a_k}\frac{n_0}{\sqrt{k}} \leq M\frac{k}{k-1} \frac{kn_0^2}{a^2_k}k^{\eta-1}\to 0, 
    \label{eq:boundsignoise2}
\end{equation}
since we are assuming that $\sqrt{k}n_0/a_k$ is upper-bounded. Therefore \eqref{eq:powersimplified} also holds in this case. In addition, we have that
\[ M_1\frac{\sqrt{k}n_0}{a_k}\leq\sqrt{k} \frac{D_k}{\sqrt{V_k}} \leq M_1\frac{\sqrt{k}n_0}{a_k}.\]
Then, as we are assuming that $M_1\leq \sqrt{k}/a_k \leq M_2$, from the monotonicity of the cumulative distribution function, it follows that $\alpha < \mathcal{P} > z_{1-\alpha}) < 1$. If additionally, \( \sqrt{k} D_k / \sqrt{V_k} \to d \in (0, \infty) \), from \eqref{eq:powersimplified} we obtain that $\mathcal{P} = \Phi(z_{\alpha} + d)$.
\item As in the previous case, \eqref{eq:powersimplified} and \eqref{eq:boundsignoise2} hold. From \eqref{eq:eq:boundssignalnoise} and since $\sqrt{k}n_0/a_k \to 0$, we get that $\sqrt{k} D_k / \sqrt{V_k} \to 0$, and therefore, from \eqref{eq:powersimplified} we get that $\mathcal{P}=\alpha.$

\end{enumerate}

\section*{Acknowledgements}

This research has been financed by research project PID2022-137818OB-I00 (Ministerio de Ciencia, Innovación y Universidades, Spain).

\bibliographystyle{apalike}
\nocite{*}
\bibliography{bibliografia.bib}

\appendix

\section{Classical Results on U-Statistics}

In this section, we present some properties of U-statistics, which will be instrumental in proving the main results later in the paper. The results and definitions not explicitly referenced can be found in Sections 5.1 and 5.2 of \cite{serfling1980approximation}.

Let $h(x_1, \ldots, x_m)$ be a symmetric kernel of degree $m$. Let $X_1,\ldots,X_n$ be a random sample with $n\geq m$. The U-statistic based on $h$ and the random sample $X_1, \ldots, X_n$ is defined as
\begin{equation}
U_n = \binom{n}{m}^{-1} \sum_{1 \leq i_1 < \cdots < i_m \leq n} h(X_{i_1}, \ldots, X_{i_m}).
    \label{eq:defustat}
\end{equation}
We say that $U_n$ is a degree-$m$ U-statistic. Let $\theta = \E\{h(X_1, \ldots, X_m)\}$. The U-statistic $U_n$ is a symmetric function of its inputs and serves as an unbiased estimator of $\theta$, i.e.,
\[
\E(U_n) = \theta.
\]

To compute the variance of $U_n$, we first define the conditional expectation of the kernel $h$ given $c$ variables
\begin{equation}
    h_c(x_1, \ldots, x_c) = \E\{h(X_1, \ldots, X_m) \mid X_1=x_1, \ldots, X_c=x_c\}.
    \label{eq:definitionhc}
\end{equation}
We define the centered version of $h_c$ as
\begin{equation}
    \widetilde{h}_c(x_1, \ldots, x_c) = h_c(x_1, \ldots, x_c) - \theta.
    \label{eq:widetildeh}
\end{equation}
Note that, $\E[h_c(X_1, \ldots, X_c)] = \theta$. The variance of both $h_c$ and $\widetilde{h}_c$ is
\begin{equation}
    \sigma_c^2 = \V\{(h_c(X_1, \ldots, X_c)\} = \V\{\widetilde{h}_c(X_1, \ldots, X_c)\}, \quad c = 1, \ldots, m.
    \label{eq:definitionsigma}
\end{equation}
 We also have that 
 for indices $i_1 < \cdots < i_r$ and $j_1 < \cdots < j_s$ sharing exactly $c$ common elements, $\#\left( \{i_1, \ldots, i_r\} \cap \{j_1, \ldots, j_s\} \right) = c$, where $\#(S)$ denotes the cardinality of $S$,
\begin{equation}
    \sigma_c^2 = \C\left\{ h_r(X_{i_1}, \ldots, X_{i_r}), h_s(X_{j_1}, \ldots, X_{j_s}) \right\}.  
    \label{eq:sigmaelementscommon}
\end{equation}
 Note that the previous equation is valid if we replace the kernels for their centered version as well. The variances $\sigma_c^2$ satisfy the ordering,
\begin{equation}
    \sigma_1^2 \leq \sigma_2^2 \leq \cdots \leq \sigma_m^2.
    \label{eq:sigmaorder}
\end{equation}
Using these variances, the variance of $U_n$ is given by
\begin{equation}
    \V(U_n) = \binom{n}{m}^{-1} \sum_{c=1}^{m} \binom{m}{c} \binom{n-m}{m-c} \sigma_c^2.
    \label{eq:varsigma}
\end{equation}

\section{Binomial coefficient properties}
\label{section:binomials}

We summarize essential properties of binomial coefficients.

\subsection*{Identity B.1 (Symmetry)}
\begin{equation}
    \binom{n}{k} = \binom{n}{n-k}
    \label{eq:bcsymettry}
\end{equation}

\subsection*{Identity B.2}

\begin{equation}
    \binom{n}{k}\binom{k}{j}=\binom{n}{j}\binom{n-j}{k-j}.
    \label{eq:nkkj}
\end{equation}

\subsection*{Identity B.3 (Vandermonde’s Identity)}
\begin{equation}
    \sum_{j=0}^{r} \binom{m}{j} \binom{n}{r-j} = \binom{m+n}{r}.
    \label{eq:vandermonde}
\end{equation}

\end{document}